\documentclass [12pt] {amsart}
\usepackage{enumerate}
\usepackage{times}
\usepackage{latexsym,amssymb, amsmath}
\usepackage{graphicx}

\usepackage{simplemargins}
\setleftmargin{2.75cm}
\setrightmargin{2.75cm}
\settopmargin{2.0in}
\setbottommargin{2.0in}
\setallmargins{2.75cm}

\def\pf{\noindent\emph{Proof: }}
\def\stop{\hfill$\Box$}

\usepackage{amsmath, amsfonts, amssymb, amsthm}
\newtheorem{thm}{Theorem}
\newtheorem{cor}[thm]{Corollary}
\newtheorem{lemma}[thm]{Lemma}

\newtheorem{defn}[thm]{Definition}
\newtheorem{prop}[thm]{Proposition}

\numberwithin{thm}{section}

\begin{document}
	
	\title[Twice Differentiable Functions with Continuous Laplacian and Unbounded Hessian]{Examples of Twice Differentiable Functions in $\mathbb{R}^n$ with Continuous Laplacian and Unbounded Hessian }
	
	\author{Yifei Pan}
	\address{Department of Mathematical Sciences\\	Purdue University Fort Wayne\\	Fort Wayne, Indiana 46805}
	\email{pan@pfw.edu}

	\author{ Yu Yan  }
	\address {Department of Mathematics and Computer Science\\ Biola University\\La Mirada, California 90639}
	\email {yu.yan@biola.edu}
	
	\begin{abstract}
	We construct examples of twice differentiable functions in $\mathbb{R}^n$ with continuous Laplacian and unbounded Hessian.	The same construction is also applicable to higher order differentiability.
	\end{abstract}

	\maketitle
	\newtheorem{Thm}{Theorem}
	\newtheorem{Lemm}{Lemma}
	\newtheorem{Cor}{Corollary}

\section{Introduction}
\label{section:introduction}

\vspace{.1in}
The standard Schauder theory states that if $\Delta u=f$ in $B_1(0) \subset \mathbb{R}^n$ and $f$ is H\"older continuous $( C^{0,\alpha}$, $ 0< \alpha <1$), then $u$ is $ C^{2,\alpha}$.  However, it fails when $\alpha=0$, that is, if $\Delta u$ is just continuous, then $u$ may not be $C^{2}$, as shown by a standard example in $\mathbb{R}^2$ (see \cite{Oton}): 
\begin{equation*}
	w(x,y)= \left\{
	\begin{array}{l l }      
		(x^2-y^2)\ln (-\ln (x^2+y^2)) & \hspace{.4in}  0<x^2+y^2\leq \frac{1}{4} , \\
		\noalign{\smallskip}
		0 & \hspace{.4in}  (x,y)=(0,0).	
	\end{array}\right.
\end{equation*}

\noindent
This function has continuous Laplacian but is not $C^2$ because it is not twice differentiable at the origin.  (Another such example can be obtained by replacing $x^2-y^2$ with $xy$.)

\vspace{.05in}
The main goal of this paper is to construct (a family of) functions that are twice differentiable everywhere with continuous Laplacian and unbounded Hessian.  These functions have only gained twice differentiablity at the origin over the above example, nevertheless, it appears that some effort is needed to achieve the gain. 

\vspace{.05in}
\begin{thm}
	\label{thm:main}
	Given any $C^2$ function $\varphi: (0, \infty) \to \mathbb{R}$ satisfying
	\begin{equation*}
		\lim _{s \to \infty} \varphi  (s) = \infty,  \hspace{.2in} \lim _{s \to \infty} \varphi ' (s) = 0,  \hspace{.2in} \lim _{s \to \infty} \varphi ''(s) = 0,
	\end{equation*}
	
	\noindent
	there is a function $u: \mathbb{R}^n \to \mathbb{R}$ depending on $\varphi$ with compact support, such that it is twice differentiable everywhere in $\mathbb{R}^n$, and it has continuous Laplacian and unbounded Hessian.  In particular, $u$ is not in $C^2(\mathbb{R}^n)$.
\end{thm}

\noindent
Obviously there are many choices of such functions $\varphi$; for example, $\varphi(s)= s^{\alpha}$ with $ 0<\alpha <1$, $\varphi(s)=\ln (s)$, or $\varphi(s)= \ln   \ln \cdots  \ln s  $ if $s>c$.

\vspace{.1in}
\noindent
As a consequence, the following is a simple application to the Dirichlet problem.

\begin{cor}
	There is a continuous function $f$ such that the unique solution of the Dirichlet problem
	$$
	\left\{
	\begin{array}{r l l}      
		\displaystyle \Delta u (x) & = & f (x)  \hspace{.4in} \text {in} \hspace{.1in} B_1(0) , \\
		\noalign{\smallskip}
		u(x) & = & 0   \hspace{.65in} \text {on} \hspace{.1in} \partial B_1(0)	
	\end{array}\right.
	$$
	is twice differentiable in $\overline{B_1(0)} $ and has unbounded Hessian. 	
\end{cor}

For any positive integer $k$, the Schauder theory also asserts that if $\Delta u$ is $C^{k, \alpha}$, then $u$ is $C^{k+2, \alpha}$. Once again, it fails when $\alpha=0$, that is, if $ \Delta u$ is just $C^{k}$, then $u $ may not be  $C^{k+2}$.  As a result of our construction we have an extension of Theorem \ref{thm:main}.

\begin{thm}
	\label{thm: kth_order}
	Given any $C^{k+2}$ function $\varphi: (0, \infty) \to \mathbb{R}$ satisfying
	\begin{equation*}
		\lim _{s \to \infty} \varphi  (s) = \infty,  \hspace{.2in} \lim _{s \to \infty} \varphi ' (s) = \cdots = \lim _{s \to \infty} \varphi ^{(k+2)}(s) = 0,
	\end{equation*}
		there is a function $u:\mathbb{R}^n \to \mathbb{R}$ depending on $\varphi$ with compact support, such that $u$ is $(k+2)$-times differentiable everywhere in $\mathbb{R}^n$, $\Delta u$ is $C^k$, but $D^{k+2}u$ is unbounded.  In particular, $u$ is not in $C^{k+2}(\mathbb{R}^n)$.
\end{thm}

\begin{cor}
	
	There is a $C^k$ function $f$ such that the unique solution of the Dirichlet problem
	$$
	\left\{
	\begin{array}{r l l}      
		\displaystyle \Delta u (x) & = & f (x)  \hspace{.4in} \text {in } \hspace{.1in} B_1(0) , \\
		\noalign{\smallskip}
		u(x) & = & 0   \hspace{.65in} \text {on} \hspace{.1in} \partial B_1(0)	
	\end{array}\right.
	$$
	is $(k+2)$-times differentiable in $\overline{B_1(0)}$, but $D^{k+2}u$ is unbounded. 	
\end{cor}

\noindent
We would like to point out a dichotomy: although the Schauder theory fails when $\alpha=0$ for each $k$, it is indeed true if $k = \infty$, since $\Delta u \in C^{\infty}$ does imply $u  \in C^{\infty}$ by the elliptic theory.

\vspace{.1in}
According to Theorem \ref{thm:main}, it would be rather natural to ask if unbounded Hessian is the only reason that hinders $u$ from being $C^2$.  Thus we propose the following problem.

\vspace{.1in}
\noindent
\textbf{Problem}:
If a function $u$ is twice differentiable everywhere, $\Delta u$ is continuous, and the Hessian of $u$ is locally bounded, then is $u$ always $C^2$?

\vspace{.05in}
\noindent
We remark that the method of construction for Theorem \ref{thm:main}  will not yield examples of twice differentiable function with continuous Laplacian and bounded Hessian without being $C^2$.  On the other hand, there are simple functions with continuous Laplacian and bounded Hessian that is twice differentiable everywhere except at the origin, such as the function (\cite{Oton-Example})
\begin{equation*}
	\phi(x, y) =	\left\{
	\begin{array}{l l }      
		(x^2-y^2)\sin (\ln (-\ln (x^2+y^2))) & \hspace{.4in}  0<x^2+y^2\leq \frac{1}{4} , \\
		\noalign{\smallskip}
		0 & \hspace{.4in}  (x,y)=(0,0).	
	\end{array}\right.
\end{equation*}

\vspace{.1in}
We also observe that if $\Delta u$ is continuous, then the modulus of continuity of $Du$ is of $o(L\log L)$. (If $\Delta u$ is just bounded, then the modulus of continuity of $Du$ is only of $O(L\log L)$.  \cite{X-Wang}) Precisely, the following is true.

\begin{prop}
	\label{thm:gradient}
	Let $u$ be a  $C^1$ solution of $\Delta u =f$, where $f$ is a continuous function on $B_1(0)$ in $\mathbb{R}^n$. Then for any $x,y \in B_{\frac{1}{2}}(0)$, 
	$$ \big | Du(x) - Du (y) \big |\leq Cd \left ( \sup_{B_1} |u| + \sup_{B_1} |f| + \int _{d}^1 \frac{\omega(r)}{r} dr \right )   , $$ where $d=|x-y|$,  $\omega(r)=\displaystyle \sup _{|x-y|<r} |f(x)-f(y)|$, and $C$ is a constant depending only on $n$.
	 
\end{prop}

\noindent
Here we notice that $$ \lim _{d \to 0} \frac{d\int _{d}^1 \frac{\omega(r)}{r} dr}{d\ln d} =0, $$which can be easily proved by considering two cases: $\displaystyle \lim _{d \to 0} \int _{d}^1 \frac{\omega(r)}{r} dr < \infty$ or  $ \displaystyle \lim _{d \to 0} \int _{d}^1 \frac{\omega(r)}{r} dr = \infty .$  It is this $o(L\log L)$ observation that motivated us to Theorem \ref{thm:main} and Theorem \ref{thm: kth_order}.


\vspace{.1in}
Theorem \ref{thm:main} and \ref{thm: kth_order} will be proved in Sections \ref{section:example_unbounded_hessian} and \ref{section:higher_order}, respectively, after a thorough study of a building block function in Section \ref{section:building_block}. One of the ideas in the construction has its origin in \cite{Coffman-Pan-Zhang} and \cite{Coffman-Pan}, where the inhomogeneous Cauchy-Riemann equation in the complex plane was considered.  Since the proof of Proposition \ref{thm:gradient} is almost identical to that for Corollary 1 in \cite{X-Wang}, we include a detailed proof in the Appendix for the convenience of the reader.

	\vspace{.2in}

\section{A Building Block Function}
\label{section:building_block}

\vspace{.1in}

In this section we will look at a function that will become a building block and provide some crucial insights for the construction of examples for Theorem \ref{thm:main}.  

\vspace{.05in}
Recall that $\varphi$ is a function satisfying
	\begin{equation}
	\label{eqn: condition_phi}
	\lim _{s \to \infty} \varphi  (s) = \infty,  \hspace{.2in} \lim _{s \to \infty} \varphi ' (s) = 0,  \hspace{.2in} \lim _{s \to \infty} \varphi ''(s) = 0.
\end{equation}

\noindent
Thus for $x=(x_1,...,x_n) \in \mathbb{R}^n$,  $\displaystyle \lim_{|x| \to 0} \varphi(-\ln |x|^2) = \infty$. Nevertheless, the product of $\varphi(-\ln |x|^2)$ with positive powers of $|x|$ is well controlled, as shown by the following two simple lemmas that will be repeatedly used later in the construction.

\begin{lemma}
	\label{lemma:x-times-varphi}
	For any  $\beta >0$ and $\varphi$ satisfying (\ref{eqn: condition_phi}),
	\begin{equation*}
		\label{eqn: z_times_varphi}
		\lim _{|x|  \to 0}  |x|^\beta  \varphi \left (-\ln  |x|^2 \right) =0.
	\end{equation*} 
	
\end{lemma}

\pf Letting $s=\frac{1}{|x|} $,  $\displaystyle   \lim _{|x| \to 0}  |x|^\beta  \varphi \left ( -\ln  |x|^2  \right )   =   \lim_{s \to \infty} \frac{\varphi (\ln s^2)}{s^{\beta}} $ is of $\frac{\infty}{\infty}$ type.  By the L'Hopital's Rule,
$$  \lim_{s \to \infty} \frac{\varphi (\ln s^2)}{s^{\beta}}   = \lim_{s \to \infty} \frac { 2\varphi ' (\ln s^2) }{\beta s^{\beta}} =0, $$ as long as $\beta >0$.

\stop

\begin{lemma}
	\label{lemma:t_times_x^t_time_varphi}
	For any $0 < \beta \leq 1$ and $\varphi$ satisfying (\ref{eqn: condition_phi}), there is a constant $B_1$ depending only on $\varphi$, such that  $$   \sup_{\substack{|x| \leq \frac{2}{3} \\ 0<\beta \leq 1}} \beta |x|^\beta \left \vert  \varphi \left (-\ln  |x|^2 \right) \right \vert \leq B_1. $$
\end{lemma}

\pf
By Lemma \ref{lemma:x-times-varphi}, we know $$ \lim _{|x|  \to 0}  \beta |x|^\beta  \varphi \left (-\ln  |x|^2 \right) =0. $$

\noindent
When $|x|=\frac{2}{3}$, $$\beta |x|^\beta \left \vert \varphi \left (-\ln  |x|^2 \right)  \right \vert = \beta \left ( \frac{2}{3} \right )^\beta  \left \vert \varphi \left (\ln \left (\frac{9}{4} \right ) \right ) \right \vert  \leq \left \vert \varphi \left (\ln \left (\frac{9}{4} \right ) \right ) \right \vert ,$$ which is a constant independent of $t$.

\vspace{.05in}
\noindent
It remains to show that the local maximum of $\beta |x|^{\beta} \left \vert \varphi \left (-\ln  |x|^2 \right)  \right \vert  $ is also bounded by a constant depending only on $\varphi$. Denote $s=\frac{1}{|x|}$, then $s > 1$, and $$\beta  |x|^{\beta} \varphi \left (-\ln  |x|^2  \right )$$ is equivalent to $$ \lambda (s) = \beta \frac{\varphi(\ln s^2) }{s^{\beta}} . $$ 

\vspace{.05in}
\noindent
We will find the local extremum of $\lambda(s)$. Because
$$\lambda '(s)= \beta \frac{2\varphi ' (\ln s^2) - \beta  \varphi(\ln s^2)}{s^{\beta +1}}, $$  a critical point 
$s_0$ must satisfy $$2\varphi ' (\ln s_0^2) = \beta  \varphi(\ln s_0 ^2).$$  At this point the local extremum of $\lambda$ is

$$
 \lambda	(s_0)  =  \beta \frac{\varphi(\ln s_0^2) }{s_0^{\beta}} = \beta \frac{\frac{2}{\beta }\varphi '(\ln s_0^2) }{s_0^{\beta}} =   \frac{2\varphi '(\ln s_0^2) }{s_0^{\beta}} ,
$$ 

\noindent
and so $$ |\lambda	(s_0) | = \left \vert   \frac{2\varphi '(\ln s_0^2) }{s_0^{\beta}} \right \vert < 2 \left \vert \varphi '(\ln s_0 ^2) \right \vert .  $$

\noindent
From (\ref{eqn: condition_phi}) we know that $\displaystyle  |\varphi '(\ln s_0^2)|  $ is bounded by a constant depending only on $\varphi$.  Thus the local maximum of $ |\lambda (s)|= \beta  |x|^{\beta} \left \vert \varphi \left (-\ln  |x|^2 \right)  \right \vert $ is also bounded by a constant depending only on $\varphi$. This completes the proof.

\stop

\vspace{.1in}
\noindent
It is worth noting that Lemma \ref{lemma:t_times_x^t_time_varphi} would not be true without the coefficient of $\beta$ in the function.   For example, let $\displaystyle \{ x_k \}$ be a sequence of points in $\mathbb{R}^n$ with $|x_k|=e^{-2^k}$, and let $\beta_k=\frac{1}{2^k}$, then $|x_k| \to 0$ and $\beta_k \to 0$ as $k \to\infty$, but $$|x_k|^{\beta_k}   \varphi \left ( -\ln  |x_k|^2  \right ) = e^{-1}   \varphi \left ( 2^{k+1}  \right ) \to \infty. $$ 
This difference will be crucial to our construction of functions with continuous Laplacian and unbounded Hessian.

\vspace{.1in}

For any $|x| \leq \frac{1}{2}$, define a function $v(x)$ by
\begin{equation}
	\label{eqn:defn_v}
	v(x) = 	\left\{
	\begin{array}{l l }      
			x_1x_2 \varphi \left (-\ln |x|^2 \right) & \hspace{.1in} 0 < |x| \leq  \frac{1}{2}, \\
			\noalign{\medskip}
				0 &   \hspace{.1in} x=0.
	\end{array}\right.
\end{equation}

\vspace{.05in}
\noindent
This function generalizes the function $w$ at the beginning of Section \ref{section:introduction} from $\mathbb{R}^2$ to $\mathbb{R}^n$.
It satisfies almost all the conditions in Theorem \ref{thm:main}, except one that it is not twice differentiable at the origin.

\begin{lemma}
	\label{lemma:v}
	The function $v$ defined by (\ref{eqn:defn_v}) has continuous Laplacian and unbounded Hessian, but it is not twice differentiable at $0$.
\end{lemma}
\vspace{.05in}
\pf
By definition, $v(x)$ is $C^2$ for all $x \neq 0$, and its derivatives are the following. (In the case $n \geq 3$, we use $i$ and $j$ to denote indices that are greater than or equal to 3.) 

\allowdisplaybreaks
\begin{eqnarray*}
	\frac{\partial v}{\partial x_1} (x) & = & x_2 \varphi \left (-\ln |x|^2 \right) - \frac{2x_1^2x_2}{|x|^2}  \varphi ' \left (-\ln |x|^2 \right), \\
	\noalign{\medskip}
	\frac{\partial v}{\partial x_2} (x) & = & x_1 \varphi \left (-\ln |x|^2 \right) - \frac{2x_1 x_2^2}{|x|^2} \varphi ' \left (-\ln |x|^2 \right), \\
	\frac{\partial v}{\partial x_i} (x) & = & - \frac {2x_1x_2x_i}{|x|^2}  \varphi ' \left (-\ln |x|^2 \right), \\
	\noalign{\medskip}
		\frac{\partial^2 v}{\partial x_1 ^2} (x) & = &  -\frac {6 x_1x_2}{|x|^2}  \varphi ' \left (-\ln |x|^2 \right)  +  \frac {4x_1^3 x_2}{|x|^4} \varphi ' \left (-\ln |x|^2 \right)  +  \frac {4x_1^3 x_2}{|x|^4} \varphi '' \left (-\ln |x|^2 \right), \\
	\noalign{\medskip}
	\frac{\partial^2 v}{\partial x_2 ^2} (x) & = &  -\frac {6x_1x_2}{|x|^2}  \varphi ' \left (-\ln |x|^2 \right)  +   \frac {4x_1 x_2^3}{|x|^4} \varphi ' \left (-\ln |x|^2 \right)   + \frac {4x_1 x_2^3}{|x|^4} \varphi '' \left (-\ln |x|^2 \right), \\
	\noalign{\medskip}
	\frac{\partial^2 v}{\partial x_i ^2} (x) & = & -\frac {2x_1x_2}{|x|^2}  \varphi '  \left (-\ln |x|^2 \right)  + \frac {4x_1x_2x_i^2}{|x|^4}  \varphi ' \left (-\ln |x|^2 \right)  +  \frac {4x_1x_2x_i^2}{|x|^4}  \varphi '' \left (-\ln |x|^2 \right), \\
	\noalign{\medskip}
	\frac{\partial^2 v}{\partial x_1 \partial x_2} (x) & = &  \varphi  \left (-\ln |x|^2 \right) - \frac{2(x_1^2+x_2^2)}{|x|^2} \varphi ' \left (-\ln |x|^2 \right) \\
	& +  & \frac {4x_1^2x_2^2}{|x|^4} \varphi ' \left (-\ln |x|^2 \right)  + \frac {4x_1^2x_2^2}{|x|^4} \varphi '' \left (-\ln |x|^2 \right ), \\
	\noalign{\medskip}
	\frac{\partial^2 v}{\partial x_1 \partial x_i} (x) & = & \frac{-2x_2x_i}{|x|^2} \varphi ' \left (-\ln |x|^2 \right) + \frac{4x_1^2x_2x_i}{|x|^4}  \varphi ' \left (-\ln |x|^2 \right)  +   \frac {4x_1^2x_2x_i}{|x|^4} \varphi '' \left (-\ln |x|^2 \right) , \\
	\noalign{\medskip}
	\frac{\partial^2 v}{\partial x_2 \partial x_i} (x) & = & \frac{-2x_1x_i}{|x|^2} \varphi ' \left (-\ln |x|^2 \right) + \frac{4x_1x_2^2x_i}{|x|^4}  \varphi ' \left (-\ln |x|^2 \right)  +   \frac {4x_1x_2^2x_i}{|x|^4} \varphi '' \left (-\ln |x|^2 \right) , \\
	\noalign{\medskip}
	\frac{\partial^2 v}{\partial x_j \partial x_i} (x) & = & -\frac{2x_1x_2}{|x|^2}\delta_{ij} \varphi ' \left (-\ln |x|^2 \right)  + \frac{4 x_1x_2x_ix_j}{|x|^4} \varphi ' \left (-\ln |x|^2 \right) + \frac{4x_1x_2x_ix_j}{|x|^4}  \varphi '' \left (-\ln |x|^2 \right) .  
\end{eqnarray*}

\vspace{.1in}
\noindent
We observe that each term in these derivatives is of the form $p(x)\varphi  \left (-\ln |x|^2 \right)$, $p(x)\varphi ' \left (-\ln |x|^2 \right)$, or $p(x)\varphi '' \left (-\ln |x|^2 \right)$, where $p(x)$ is homogeneous in $x$. For first derivatives, the degree of homogeneity is 1, and for second derivatives, the degree of homogeneity is 0. Because of this, by Lemma \ref{lemma:x-times-varphi} and the choice of $\varphi$, $$ \displaystyle \lim _{|x| \to 0} \frac{\partial v}{\partial x_1} (x) =  \lim _{|x| \to 0} \frac{\partial v}{\partial x_2} (x) =  \lim _{|x| \to 0} \frac{\partial v}{\partial x_i} (x)  =0,$$ 
	
	\begin{eqnarray*} 
		& & \lim _{|x| \to 0} \frac{\partial ^2 v}{\partial x_1^2} (x) =  \lim _{|x| \to 0} \frac{\partial ^2 v}{\partial x_2^2} (x) =  \lim _{|x| \to 0} \frac{\partial ^2 v}{\partial x_i^2} (x)  \\
			\noalign{\medskip}
		& = & \lim _{|x| \to 0} 	\frac{\partial^2 v}{\partial x_1 x_i} (x) = \lim _{|x| \to 0} 	\frac{\partial^2 v}{\partial x_2 x_i} (x)= \lim _{|x| \to 0} 	\frac{\partial^2 v}{\partial x_i x_j} (x) = 0.
	\end{eqnarray*} 

\vspace{.05in}
	\noindent
Thus all the first and second derivatives of $v$ approach $0$ as $|x| \to 0$, except for  
\begin{eqnarray*}
	\label{eqn:d^2v_dx_1_dx_2}
	\frac{\partial^2 v}{\partial x_1 \partial x_2} (x) & = &  \varphi  \left (-\ln |x|^2 \right) - \frac{2(x_1^2+x_2^2)}{|x|^2} \varphi ' \left (-\ln |x|^2 \right) \\
& +  & \frac {4x_1^2x_2^2}{|x|^4} \varphi ' \left (-\ln |x|^2 \right)  + \frac {4x_1^2x_2^2}{|x|^4} \varphi '' \left (-\ln |x|^2 \right ).
\end{eqnarray*}

\vspace{.05in}
\noindent
As $|x| \to 0$, its first term goes to $\infty$ and all the other terms go to 0, thus $ \frac{\partial^2 v}{\partial x_1 \partial x_2} (x)$ is unbounded near the origin, which causes the Hessian of $v$ to be unbounded.  

\vspace{.05in}
\noindent
On the other hand, because each diagonal entry of the Hessian has a removable discontinuity at the origin,  
\begin{eqnarray*}
	\Delta v(x) & = & 	\frac{\partial^2 v}{\partial x_1 ^2} (x) + 	\frac{\partial^2 v}{\partial x_2 ^2} (x) + \sum _{i=3}^n 	\frac{\partial^2 v}{\partial x_i ^2} (x) \\
	& = & - \frac {( 2n+8 ) x_1x_2}{|x|^2}  \varphi ' \left (-\ln |x|^2 \right)  +   \frac {4x_1 x_2 }{|x|^2} \Big ( \varphi ' \left (-\ln |x|^2 \right)  +\varphi '' \left (-\ln |x|^2 \right) \Big ) \\  
		\noalign{\medskip}
	& \to & 0 \hspace{.2in} \text{as} \,\, |x| \to 0.
\end{eqnarray*}

\vspace{.05in}
\noindent
Lastly, we check the differentiablity of $v$ at 0. It is differentiable because by Lemma \ref{lemma:x-times-varphi},
$$
\frac{\left \vert v(x)-v(0)  \right \vert}{|x|}  =  	\frac{\left \vert 	x_1x_2 \varphi \left (-\ln |x|^2 \right)  \right \vert}{|x|} \leq |x| \left \vert    \varphi \left (-\ln |x|^2 \right) \right \vert \to 0 \hspace{.2in} \text{as} \,\, |x| \to 0.
$$

\vspace{.05in}
\noindent
Therefore, all first derivatives of $v$ equal 0 at the origin, and $v$ is $C^1$ throughout $\mathbb{R}^n$. 
Computing the partial derivative $\frac{\partial^2 v}{\partial x_1 ^2} (0)$ by definition, we have
\begin{equation*}
	\frac{\partial^2 v}{\partial x_1 ^2} (0) = \lim_{h \to 0} \dfrac{\frac{\partial v}{\partial x_1 } ((h,0,...,0)) - \frac{\partial v}{\partial x_1 } ((0,...,0))  }{h} =  \lim_{h \to 0} \frac{0-0}{h} =0.
\end{equation*}
Similarly, we also have $$\frac{\partial^2 v}{\partial x_2 ^2} (0) = \cdots = \frac{\partial^2 v}{\partial x_n ^2} (0)=0.$$

\vspace{.05in}
\noindent
Therefore, $\Delta v (0)=0$.  Consequently, $\Delta v$ is continuous at $0$.
 However, $v$ is not twice differentiable at 0.  To see that, we check the differentiability of $	\frac{\partial v}{\partial x_1 }$ at $0$:
\begin{equation*}
	\label{eqn:dv_dx_1-difference_quotient}
	\frac{\left \vert 	\frac{\partial v}{\partial x_1 } (x)-	\frac{\partial v}{\partial x_1 } (0)  \right \vert}{|x|}  =  	\frac{\left \vert x_2 \varphi \left (-\ln |x|^2 \right) - \frac{2x_1^2x_2}{|x|^2}  \varphi ' \left (-\ln |x|^2 \right )\right \vert}{|x|} .
\end{equation*}

\vspace{.05in}
\noindent
Let $x_1=x_3=\cdots=x_n=0$, then $|x_2|=|x|$ and
$$
\frac{\left \vert 	\frac{\partial v}{\partial x_1 } (x)-	\frac{\partial v}{\partial x_1 } (0)  \right \vert}{|x|}  =  	\frac{\left \vert x_2 \varphi \left (-\ln |x_2|^2 \right) \right \vert}{|x_2|}= \left \vert  \varphi \left (-\ln |x_2|^2 \right) \right \vert \to \infty \hspace{.2in} \text{as} \,\, x_2 \to 0,
$$

\vspace{.05in}
\noindent
hence $ 	\frac{\partial v}{\partial x_1 } $ is not differentiable at $0$.  Similarly, $	\frac{\partial v}{\partial x_2 } $ is not differentiable at $0$ either.

\vspace{.05in}
\noindent
Interestingly,
$$ \frac{\left \vert 	\frac{\partial v}{\partial x_i } (x)-	\frac{\partial v}{\partial x_i } (0)  \right \vert}{|x|}  =  	\frac{\left \vert - \frac {2x_1x_2x_i}{|x|^2}  \varphi ' \left (-\ln |x|^2 \right)\right \vert}{|x|} \leq 2\left \vert \varphi ' \left (-\ln |x|^2 \right)\right \vert \to 0 \hspace{.2in} \text{as} \,\, x \to 0,$$
thus 
$ 	\frac{\partial v}{\partial x_i } $ for all $i \geq 3$ are differentiable at $0$.

\vspace{.05in}
\noindent
Therefore, $v$ fails to be twice differentiable at $0$ because $	\frac{\partial v}{\partial x_1 } $ and $	\frac{\partial v}{\partial x_2 } $ are not differentiable at $0$.
This completes the proof of Lemma \ref{lemma:v}.

\stop




\vspace{.2in}

\section{Construction for Theorem \ref{thm:main}}
\label{section:example_unbounded_hessian}

\vspace{.1in}

In this section, we will first ``smooth out" $v$ into a function that is $C^2$ at the origin, then we will combine a sequence of such functions through scaling and translation to create a desired function that is twice differentiable everywhere with continuous Laplacian and unbounded Hessian, thus proving Theorem \ref{thm:main}.


\begin{defn}
	\label{defn: u_t}
	\vspace{.1in}
	\noindent 
	Let $\eta: [0, \infty) \to [0,1]$ be a fixed, non-increasing $C^{\infty}$ function such that
	\begin{equation}
		\label{eqn: condition_eta}
		\eta(s) \equiv 1 \,\, \text {for } \,\, 0 \leq s \leq \frac{1}{2}  \hspace{.2in}   \text{and} \hspace{.2in} \eta(s) \equiv 0 \,\, \text {for } \,\, s \geq \frac{2}{3}  . 
	\end{equation}
	
\vspace{.05in}
	\noindent 
	For any $0<t\leq \frac{1}{2}$, define a function $u_t: \mathbb{R}^n \to \mathbb{R} \,\, (n \geq 2)$  by
	\begin{equation}
		\label{eqn: defn_u_t}
		u_t(x) = 	\left\{
		\begin{array}{l l l}      
			0 &   \hspace{.1in} x=0, \\
			\noalign{\medskip}
			\eta(|x|)x_1 x_2|x|^{2t}\varphi  (-\ln  |x|^2  ) & \hspace{.1in} 0 < |x|< 1, \\
			\noalign{\medskip}
			0 &	\hspace{.1in} |x| \geq 1.
		\end{array}\right.
	\end{equation}
	
\end{defn}

\noindent 
It follows immediately from Lemma \ref{lemma:x-times-varphi} and (\ref{eqn: condition_eta}) that $u_t$ is continuous everywhere. Actually, it can be shown that $u_t \in C^2 \left ( \mathbb{R}^n  \right )$, but we will not verify it here because it is not to be used in our construction. 



\vspace{.05in}
What will be essential to our construction is the fact that all the first derivatives of $u_t$ and second derivatives of the form $ \frac{\partial ^2 u_t}{\partial x_j ^2}  $ are uniformly bounded by constants independent of $t$. However,  that is not the case for $ \frac{\partial ^2 u_t}{\partial x_1 \partial x_2} $, as will be shown later in this section.  

\begin{lemma}
	\label{lemma:u_t_derivative_bounded} 
	There are constants $B_2$ and $B_3$ depending only on $\eta$ and $\varphi$, such that 
	
	\begin{equation}
		\label{eqn:u_t_derivative_bounded}
		\sup _{x \in \mathbb{R}^n} \left \vert \frac{\partial  u_t}{\partial x_j} (x) \right \vert \leq B_2 \hspace{.4in} \text{for} \hspace{.2in} j=1, ..., n.
	\end{equation}

	\begin{equation}
	\label{eqn:u_t_derivative_bounded_second_diagonal}
	\sup _{x \in \mathbb{R}^n} \left \vert \frac{\partial ^2 u_t}{\partial x_j^2} (x) \right \vert \leq B_3 \hspace{.4in} \text{for} \hspace{.2in} j=1, ..., n.
\end{equation}

\end{lemma}

\vspace{.05in}
\pf By the definition of $u_t$, $$ \sup _{x \in \mathbb{R}^n} \left \vert \frac{\partial u_t}{\partial x_j} (x) \right \vert =  \sup _{|x| \leq 1} \left \vert \frac{\partial u_t}{\partial x_j} (x) \right \vert \hspace{.4in} \text{and} \hspace{.4in} \sup _{x \in \mathbb{R}^n} \left \vert \frac{\partial ^2 u_t}{ \partial x_j^2 } (x) \right \vert  =  \sup _{|x| \leq 1} \left \vert \frac{\partial ^2 u_t}{\partial x_j^2} (x) \right \vert. $$

\vspace{.1in}
\noindent
First, we need to compute the derivatives.  When $0< |x| < 1$,

 \begin{eqnarray}
	\label{eqn:derivatives_u_t_x_1}
	\frac{\partial u_t}{\partial x_1} (x) & = &   \eta ' (|x|) x_1^2x_2|x|^{2t-1} \varphi  (-\ln  |x|^2  ) +  \eta  (|x|) x_2|x|^{2t}  \varphi  (-\ln  |x|^2  )  \\
	&  + & (2t) \eta (|x|)x_1^2 x_2 |x|^{2t-2} \varphi  (-\ln  |x|^2  )   - 2\eta (|x|)x_1 ^2 x_2|x|^{2t-2} \varphi '  (-\ln  |x|^2  ).  \nonumber
\end{eqnarray}

\begin{eqnarray}
	\label{eqn:derivatives_u_t_x_2}
	\frac{\partial u_t}{\partial x_2} (x) & = &   \eta ' (|x|) x_1 x_2^2|x|^{2t-1}\varphi  (-\ln  |x|^2  ) +  \eta  (|x|) x_1|x|^{2t} \varphi  (-\ln x|^2  )  \\
	&  + & (2t) \eta (|x|)x_1 x_2^2 |x|^{2t-2}\varphi   (-\ln  |x|^2  )   - 2\eta (|x|)x_1  x_2^2 |x|^{2t-2}\varphi '  (-\ln |x|^2  ).  \nonumber
\end{eqnarray}

\vspace{.05in}
\noindent
In the case $n \geq 3$, for any $ i \geq 3$,
\begin{eqnarray}
	\label{eqn:derivatives_u_t_x_i}
	\frac{\partial u_t}{\partial x_i} (x) & = &   \eta ' (|x|) x_1 x_2 x_i|x|^{2t-1}\varphi  (-\ln  |x|^2  )  + (2t) \eta (|x|)x_1 x_2 x_i |x|^{2t-2}\varphi   (-\ln  |x|^2  )   \\
	& - & 2\eta (|x|)x_1  x_2 x_i |x|^{2t-2}\varphi '  (-\ln  |x|^2  ).  \nonumber 
\end{eqnarray}

\vspace{.05in}
\noindent
The first term in (\ref{eqn:derivatives_u_t_x_1}) is bounded by $$C|x|^{2t+2} \left \vert \varphi   (-\ln |x|^2  ) \right \vert, $$ where $C$ depends only on $\eta$.  Because $|x| \leq 1$, we have $ |x|^{2t+2} \leq |x|^2,$ so 
$$ |x|^{2t+2} \left \vert \varphi   (-\ln |x|^2  ) \right \vert \leq |x|^{2} \left \vert \varphi   (-\ln |x|^2  ) \right \vert.  $$  By Lemma \ref{lemma:x-times-varphi}, $|x|^2 \varphi   (-\ln |x|^2  )$ has a removable discontinuity at $0$, therefore on the closed set $|x| \leq 1$ it is bounded by a constant depending only on $\varphi$.  Thus the first term in (\ref{eqn:derivatives_u_t_x_1}) is bounded by a constant depending only on $\eta$ and $\varphi$.   Similarly, we can prove that the second and third terms in  (\ref{eqn:derivatives_u_t_x_1}) are also bounded by constants depending only on $\eta$ and $\varphi$. The last term, $$2\eta (|x|)x_1 ^2 x_2|x|^{2t-2} \varphi '  (-\ln  |x|^2  ), $$ is bounded by  $$C|x|^{2t+1} \left \vert \varphi '  (-\ln |x|^2  ) \right \vert,$$ where $C$ depends only on $\eta$.  It is further bounded by $$C|x| \left \vert \varphi '  (-\ln |x|^2  ) \right \vert$$ since $|x|\leq 1$.  Because of (\ref{eqn: condition_phi}), we know $ |x| \varphi '  (-\ln |x|^2  ) $ has a removable discontinuity at $0$, therefore on the closed set $|x| \leq 1$ it is bounded by a constant depending only on $\varphi$.  Thus the last term in (\ref{eqn:derivatives_u_t_x_1}) is also bounded by a constant depending only on $\eta$ and $\varphi$.  
   Therefore,  $	\frac{\partial u_t}{\partial x_1}$ is bounded by a constant depending on $\eta$ and $\varphi$ only.  In the same way, we can prove that $	\frac{\partial u_t}{\partial x_2}$ and $\frac{\partial u_t}{\partial x_i} (i \geq 3)$ are also bounded by a constant depending on $\eta$ and $\varphi$ only. This proves (\ref{eqn:u_t_derivative_bounded}).

\vspace{.1in}
\noindent
Next we prove (\ref{eqn:u_t_derivative_bounded_second_diagonal}).
\allowdisplaybreaks
\begin{eqnarray}
	\label{eqn:2nd_derivatives_u_t_x_1}
	\frac{\partial ^2 u_t}{\partial x_1^2} (x) 
	& = &   \eta '' (|x|)   x_1^3 x_2 |x|^{2t-2} \varphi  (-\ln  |x|^2   ) +   2 \eta ' (|x|)  x_1 x_2 |x|^{2t-1} \varphi  (-\ln  |x|^2   ) \nonumber \\
	& + & (4t-1)  \eta ' (|x|)   x_1^3 x_2 |x|^{2t-3} \varphi  (-\ln  |x|^2   )
	 +   \eta ' (|x|)   x_1 x_2 |x|^{2t-1} \varphi  (-\ln  |x|^2  ) \nonumber \\
	 &	+ & 6t \eta  (|x|)   x_1 x_2 |x|^{2t-2} \varphi  (-\ln  |x|^2  )  -  2 \eta  (|x|)   x_1 x_2 |x|^{2t-2} \varphi '  (-\ln |x|^2  )  \\ 
	 & + & 2t(2t-2) \eta  (|x|)   x_1^3 x_2 |x|^{2t-4} \varphi  (-\ln |x|^2  )  -   4t \eta  (|x|)   x_1 ^3 x_2 |x|^{2t-4} \varphi '  (-\ln ( |x|^2  ) \nonumber \\
	 & - &  4 \eta '  (|x|)   x_1^3 x_2 |x|^{2t-3} \varphi '  (-\ln  |x|^2  )  -  6 \eta  (|x|)   x_1 x_2 |x|^{2t-2} \varphi '  (-\ln  |x|^2  ) \nonumber \\
	 & - &  2(2t-2)  \eta (|x|)   x_1^3 x_2 |x|^{2t-4} \varphi '  (-\ln |x|^2  )  +  4 \eta  (|x|)   x_1^3 x_2 |x|^{2t-4} \varphi ''  (-\ln  |x|^2  ). \nonumber
\end{eqnarray}
\allowdisplaybreaks
\begin{eqnarray}
	\label{eqn:2nd_derivatives_u_t_x_2}
	\frac{\partial ^2 u_t}{\partial x_2^2} (x) & = &   \eta '' (|x|)   x_1 x_2^3 |x|^{2t-2} \varphi  (-\ln  |x|^2   )   + 2 \eta ' (|x|)  x_1 x_2 |x|^{2t-1} \varphi  (-\ln  |x|^2   ) \nonumber \\
	& + & (4t-1)  \eta ' (|x|)   x_1 x_2 ^3|x|^{2t-3} \varphi  (-\ln  |x|^2   )
 +  \eta ' (|x|)   x_1 x_2 |x|^{2t-1} \varphi  (-\ln  |x|^2  ) \nonumber \\
 &	+ &  6t \eta  (|x|)   x_1 x_2 |x|^{2t-2} \varphi  (-\ln  |x|^2  )   -  2 \eta  (|x|)   x_1 x_2 |x|^{2t-2} \varphi '  (-\ln  |x|^2  ) \nonumber \\
 & + &  2t(2t-2) \eta  (|x|)   x_1 x_2^3 |x|^{2t-4} \varphi  (-\ln  |x|^2  )  -   4t \eta  (|x|)   x_1  x_2 ^3 |x|^{2t-4} \varphi '  (-\ln  |x|^2  ) \\
 & - &  4 \eta '  (|x|)   x_1 x_2 ^3 |x|^{2t-3} \varphi '  (-\ln  |x|^2  )  -  6 \eta  (|x|)   x_1 x_2 |x|^{2t-2} \varphi '  (-\ln  |x|^2  ) \nonumber \\
 & - & 2(2t-2)  \eta (|x|)   x_1 x_2^3 |x|^{2t-4} \varphi '  (-\ln  |x|^2  )   +  4 \eta  (|x|)   x_1 x_2 ^3 |x|^{2t-4} \varphi ''  (-\ln  |x|^2  ).  \nonumber
\end{eqnarray}

\vspace{.05in}
\noindent
In the case $n \geq 3$, for any $ i \geq 3$,
\allowdisplaybreaks
\begin{eqnarray}
	\label{eqn:2nd_derivatives_u_t_x_i}
	\frac{\partial ^2 u_t}{\partial x_i^2} (x) & = &   \eta '' (|x|)   x_1 x_2 x_i^2 |x|^{2t-2} \varphi  (-\ln  |x|^2   )   + \eta ' (|x|)  x_1 x_2 |x|^{2t-1} \varphi  (-\ln  |x|^2   ) \nonumber  \\
	& + & (4t-1)  \eta ' (|x|)   x_1 x_2 x_i^2|x|^{2t-3} \varphi  (-\ln  |x|^2   )
	-2  \eta ' (|x|)   x_1 x_2 x_i^2 |x|^{2t-3} \varphi '  (-\ln  |x|^2  )  \nonumber \\
	& +  & 	 2t \eta  (|x|)   x_1 x_2 |x|^{2t-2} \varphi (-\ln  |x|^2 )  +  2t(2t-2) \eta  (|x|)   x_1 x_2 x_i^2 |x|^{2t-4} \varphi (-\ln  |x|^2  ) \\ 
	& + & (4-8t) \eta  (|x|)   x_1 x_2x_i^2 |x|^{2t-4} \varphi ' (-\ln  |x|^2  )  -     2 \eta '  (|x|)   x_1 x_2 x_i^2 |x|^{2t-3} \varphi '  (-\ln  |x|^2 ) \nonumber \\
	&  - & 2 \eta  (|x|)   x_1 x_2 |x|^{2t-2} \varphi ' (-\ln  |x|^2  )  +  4 \eta  (|x|)   x_1 x_2 x_i^2 |x|^{2t-4} \varphi ''  (-\ln |x|^2  ). \nonumber
\end{eqnarray}

\vspace{.05in}
\noindent
To prove (\ref{eqn:u_t_derivative_bounded_second_diagonal}) we need to estimate each term of (\ref{eqn:2nd_derivatives_u_t_x_1}), (\ref{eqn:2nd_derivatives_u_t_x_2}), and (\ref{eqn:2nd_derivatives_u_t_x_i}).  

\vspace{.1in}
\noindent
We start with (\ref{eqn:2nd_derivatives_u_t_x_1}).  Note that $|x| \leq 1$, the first term is bounded by $$C|x|^2  \left \vert \varphi   (-\ln |x|^2  ) \right \vert,$$ the 2nd through 4
th terms are bounded by $$C|x| \left \vert \varphi   (-\ln |x|^2  ) \right \vert ,$$  the 9th term is bounded by $$C|x|\left \vert \varphi '  (-\ln |x|^2  ) \right \vert, $$ and the 12th term is bounded by $$C \left \vert \varphi ''  (-\ln |x|^2  ) \right \vert ,$$ where $C$ is a constant depending only on $\eta$.  By the same argument as that in the above proof for the first derivatives, all these functions are bounded uniformly by a constant depending only on $\eta$ and $\varphi$.

\vspace{.1in}
\noindent
The 6th, 8th, 10th and 11th terms are bounded by $$C|x|^{2t} \left \vert\varphi '  (-\ln |x|^2  ) \right \vert,$$ which is further bounded by  $$ C \left \vert\varphi '  (-\ln |x|^2 ) \right \vert. $$ 

\noindent
By (\ref{eqn: condition_phi}), $\varphi '   (-\ln |x|^2  )$ has a removable discontinuity at $0$, therefore on the closed set $|x| \leq 1$ it is bounded by a constant depending only on $\varphi$. 

\vspace{.1in}
\noindent
The 5th term, $$6t \eta  (|x|)   x_1 x_2 |x|^{2t-2} \varphi  (-\ln  |x|^2  ),$$ and the 7th term, $$ 2t(2t-2) \eta  (|x|)   x_1^3 x_2 |x|^{2t-4} \varphi  (-\ln |x|^2  ),$$ are bounded by $$Ct|x|^{2t}\left \vert  \varphi  (-\ln |x|^2  ) \right \vert.  $$ By Lemma \ref{lemma:t_times_x^t_time_varphi}, $$  2t|x|^{2t}\left \vert  \varphi  (-\ln |x|^2  ) \right \vert \leq  B_1, $$ where $B_1$ depends only on $\varphi$. Hence the 5th and 7th terms are bounded by a constant depending on $\eta$ and $\varphi$.  Therefore, we have proved that all the terms in (\ref{eqn:2nd_derivatives_u_t_x_1}) are uniformly bounded by a constant independent of $t$.  

\vspace{.1in}
\noindent
All the terms in (\ref{eqn:2nd_derivatives_u_t_x_2}) and (\ref{eqn:2nd_derivatives_u_t_x_i}) can be estimated in the same way, so this completes the proof of (\ref{eqn:u_t_derivative_bounded_second_diagonal}).

\stop

\vspace{.1in}

 Now we are ready to construct the main function, $u$, by ``piecing together" a sequence of functions $u_{t_k}$ as follows.

\vspace{.05in}
\noindent
Choose two decreasing sequences of numbers $R_k \to 0$ and $r_k \to 0$, such that $$ R_k > r_k, $$ and for geometric reasons that will be explained later we also require
\begin{equation}
	\label{eqn:R_k_r_k}
	R_k -r_k > R_{k+1}+r_{k+1};
\end{equation}
for example, we may choose $R_k=10^{-k}$ and $r_k=10^{-(k+1)}$.

\vspace{.1in}
\noindent
We use $\zeta_0$ to denote the point $\left ( \frac{1}{\sqrt{2}}, \frac{1}{\sqrt{2}}, ..., \frac{1}{\sqrt{2}}  \right ) $ in $\mathbb{R}^n$ and choose a sequence $\{t_k\}$ such that  $0< t_k < \frac{1}{4}$ and $\displaystyle \lim _{k \to \infty} t_k =0$.  Define the function $u(x)$ by  
\begin{equation}
	\label{defn:u}
	u(x) = \sum _{k=1}^{\infty} \epsilon_k r_k ^2 u_{t_k} \left ( \frac{x-R_k \zeta_0}{r_k} \right ) ,
\end{equation}
where the only conditions on $\epsilon_k$ for now are $\epsilon_k>0  $ and $\displaystyle \lim _{k \to \infty} \epsilon_k =0$,  we do not need to assign specific values to $\epsilon_k$ until near the end of this section.

\vspace{.05in}
\noindent
Condition (\ref{eqn:R_k_r_k}) ensures that the balls centered at the points $R_k\zeta_0$ with radii $r_k$ are mutually disjoint.  For each $k \in \mathbb{N}$, let $B_k$ be the ball
centered at the point $R_k \zeta_0$ with radius $ \frac{2}{3}r_k$, then these $B_k$ are also mutually disjoint.  By (\ref{eqn: condition_eta}) and (\ref{eqn: defn_u_t}), the support of each function $u_{t_k} \left ( \frac{x-R_k \zeta_0}{r_k} \right )$ is the ball $\{x \in \mathbb{R}^n: |x-R_k \zeta_0| \leq \frac{2}{3}r_k \}$, which is $B_k$.
Therefore, although the definition of $u(x)$ appears to be an infinite sum, it actually is only a single term.   For any given $x \in \mathbb{R}^n$, if $x$ is not in any of the $B_k$, then $$u(x)=0,$$ otherwise $$u(x)=\epsilon_k r_k ^2 u_{t_k} \left ( \frac{x-R_k \zeta_0}{r_k} \right ) \hspace{.2in} \text{for some} \,\, k.$$

\vspace{.05in}
\noindent
As $k \to \infty$, the radius of $B_k$ goes down to 0 and its center moves toward the origin, but none of the balls $B_k$ contains the origin. In fact, for any $j=1,.., n$, the $x_j$-th coordinate hyperplane does not intersect any of the ball $B_k$.  To see this, let $$(x_1,...,x_{j-1}, 0, x_{j+1},...x_n)$$ be an arbitrary point on  the $x_j$-th coordinate hyperplane. The distance from this point to the center of the ball, $R_k\zeta_0 = \left ( \frac{R_k}{\sqrt{2}}, ..., \frac{R_k}{\sqrt{2}}  \right ) $, is 

$$ \sqrt{ \left( \frac{R_k}{\sqrt{2}} -x_1 \right )^2+ \cdots + \left( \frac{R_k}{\sqrt{2}} -x_{j-1 }\right )^2+ \left( \frac{R_k}{\sqrt{2}} \right )^2 + \left( \frac{R_k}{\sqrt{2}} -x_{j+1} \right )^2 + \cdots + \left( \frac{R_k}{\sqrt{2}} -x_n \right )^2 } $$
$$ \geq   \frac{R_k}{\sqrt{2}}> \frac{r_k}{\sqrt{2}} > \frac{2}{3}r_k, \hspace{6in}
$$

\vspace{.05in}
\noindent
since $\frac{1}{\sqrt{2}} \approx 0.71$ and $ \frac{2}{3} \approx 0.67$.  Thus $u = 0$ on all of the $n$ coordinate hyperplanes, and consequently $u(0)=0.$ By an argument similar to that in the proof of Lemma \ref{lemma:u_t_derivative_bounded}, we can show that $u_{t_k}$ is uniformly bounded by a constant independent of $t_k$, hence by definition $\displaystyle \lim _{|x| \to 0} u(x)=0$.  Therefore $u$ is continuous at the origin, and thus continuous everywhere in $\mathbb{R}^n$.

\vspace{.05in}
\begin{lemma}
	\label{lemma:u_differentiable}
	The function $u(x)$ as defined in (\ref{defn:u}) is twice differentiable everywhere in $\mathbb{R}^n$, and all its first and second order partial derivatives at the origin are equal to 0. 
\end{lemma}

\pf 
By definition $u(x)$ is $C^2$ for all $x \neq 0$, so we only need to show it is twice differentiable at the origin.

\vspace{.1in}
\noindent
Because $u=0$ on all the coordinate hyperplanes,
\begin{equation*}
	\label{eqn:u_derivative_0}
	\frac{\partial u}{\partial x_j} (0)=0.
\end{equation*}

\vspace{.05in}
\noindent
From (\ref{defn:u}), we have 
\begin{equation}
	\label{eqn:derivative_u}
	\frac{\partial u}{\partial x_j} (x)= \sum _{k=1}^{\infty} \epsilon_k r_k  \frac{\partial u_{t_k}}{\partial x_j} \left ( \frac{x-R_k \zeta_0}{r_k} \right ). \nonumber
\end{equation}

\vspace{.05in}
\noindent
Recall that the balls $B_k$ are mutually disjoint, so for any given $x \in \mathbb{R}^n$, either $$\frac{\partial u}{\partial x_j}(x)=0,$$ or $$\frac{\partial u}{\partial x_j}(x) = \epsilon_k r_k  \frac{\partial u_{t_k}}{\partial x_j} \left ( \frac{x-R_k \zeta_0}{r_k} \right ) \hspace{.2in} \text{for some} \,\, k.$$ 

\vspace{.05in}
\noindent
By (\ref{eqn:u_t_derivative_bounded}) in Lemma \ref{lemma:u_t_derivative_bounded}, $$ \left \vert \epsilon_k r_k  \frac{\partial u_{t_k}}{\partial x_j} \left ( \frac{x-R_k \zeta_0}{r_k} \right )  \right \vert  \leq \epsilon_k r_k  B_2 \to 0 \hspace{.2in} \text{as} \,\, k \to \infty. $$ 

\vspace{.05in}
\noindent
Therefore, $$ \lim _{|x| \to 0}  \frac{\partial u}{\partial x_j}(x) =0, $$ which implies $u$ is $C^1$ at the origin.  By definition it is $C^1$ for all $x \neq 0$, thus $u$ is $C^1$ everywhere in $\mathbb{R}^n$.
Furthermore, for any $j=1,..,n$,
$$ 	\frac{\left \vert \frac{\partial u}{\partial x_j}(x)-\frac{\partial u}{\partial x_j}(0) \right \vert }{|x|} = 	\frac{\left \vert \frac{\partial u}{\partial x_j}(x) \right \vert }{|x|}. $$

\vspace{.05in}
\noindent
If $\displaystyle \frac{\partial u}{\partial x_j}(x)=0$, then $\displaystyle  \frac{\left \vert \frac{\partial u}{\partial x_j}(x) \right \vert }{|x|} =0$.
If $ \displaystyle \frac{\partial u}{\partial x_j}(x) = \epsilon_k r_k  \frac{\partial u_{t_k}}{\partial x_j} \left ( \frac{x-R_k \zeta_0}{r_k} \right )$ for some $k$, then by (\ref{eqn:u_t_derivative_bounded})

$$
	\frac{\left \vert \frac{\partial u}{\partial x_j}(x) \right \vert }{|x|}  \leq  \frac{\epsilon_k r_k B_2 }{R_k - \frac{2}{3}r_k} = \frac{\epsilon_k  B_2 }{\frac{R_k}{r_k} - \frac{2}{3}} <  \frac{\epsilon_k  B_2 }{1 - \frac{2}{3} } \to 0  \hspace{.2in} \text{as} \,\, k \to \infty. 
$$

\vspace{.1in}
\noindent
Hence we know that $$ \lim _{|x| \to 0} \frac{\left \vert \frac{\partial u}{\partial x_j}(x)-\frac{\partial u}{\partial x_j}(0) \right \vert }{|x|}  =0, $$ 

\vspace{.1in}
\noindent
which means $u$ is twice differentiable at the origin, where all its second order partial derivatives are 0. This completes the proof.

\stop

\vspace{.05in}
Lastly, we will show that $u$ has continuous Laplacian but unbounded Hessian.

\vspace{.05in}
\begin{lemma}
	\label{lemma:u_not_C^2_Laplacian continuous}
	The function $u(x)$ as defined in (\ref{defn:u}) has continuous Laplacian everywhere in $\mathbb{R}^n$, and the partial derivative $\frac{\partial^2 u}{\partial x_1 \partial x_2}$ is unbounded near the origin.
\end{lemma}

\pf
Because $u=0$ on all the coordinate hyperplanes,
\begin{equation*}
	\label{eqn:u_derivative_0}
	\frac{\partial ^2 u}{\partial x_j^2} (0)=0.
\end{equation*}

\noindent
 For any given $x \in \mathbb{R}^n$, either $$\frac{\partial^2 u}{\partial x_j^2}(x)=0,$$ or $$\frac{\partial ^2 u}{\partial x_j^2}(x) = \epsilon_k   \frac{\partial^2 u_{t_k}}{\partial x_j^2} \left ( \frac{x-R_k \zeta_0}{r_k} \right ) \hspace{.2in} \text{for some} \,\, k.$$ 
 By (\ref{eqn:u_t_derivative_bounded_second_diagonal}) in Lemma \ref{lemma:u_t_derivative_bounded}, $$   \epsilon _k \left \vert  \frac{\partial^2 u_{t_k}}{\partial x_j^2} \left ( \frac{x-R_k \zeta_0}{r_k} \right )   \right \vert  \leq \epsilon _k  B_3 \to 0 \hspace{.2in} \text{as} \,\, k \to \infty. $$ Thus $\displaystyle \lim_{|x| \to 0} \frac{\partial ^2 u}{\partial x_j^2}(x) = 0 $, which implies that $\displaystyle \frac{\partial^2 u}{\partial x_j^2}(x)$ is continuous at 0.  Since it is also continuous for all $x \neq 0$, it is continuous everywhere.
This proves that $\displaystyle \Delta u = \sum _{j=1}^n \frac{\partial^2 u}{\partial x_j^2}$ is continuous everywhere in $\mathbb{R}^n$.

\vspace{.05in}
\noindent
Next we will show $\displaystyle \frac{\partial^2 u}{\partial x_1 \partial x_2}$ is unbounded near the origin.
For general $t$,
\allowdisplaybreaks
\begin{eqnarray}
	\label{eqn:2nd_derivatives_u_t_x_1_x_2}
	\frac{\partial ^2 u_t}{\partial x_1 \partial x_2} (x) & = &   \eta '' (|x|)   x_1^2 x_2 ^2 |x|^{2t-2} \varphi  (-\ln  |x|^2   )   + \eta ' (|x|) \left ( x_1^2 + x_2^2 \right )  |x|^{2t-1} \varphi  (-\ln  |x|^2   ) \nonumber  \\
	& + & (4t-1)  \eta ' (|x|)   x_1^2 x_2^2|x|^{2t-3} \varphi  (-\ln  |x|^2   )
	- 4 \eta ' (|x|)   x_1^2 x_2^2 |x|^{2t-3} \varphi '  (-\ln  |x|^2  )  \nonumber \\
	& +  & 	 \eta  (|x|)    |x|^{2t} \varphi (-\ln  |x|^2 )  +  2t \eta  (|x|) \left ( x_1^2 + x_2^2 \right )    |x|^{2t-2} \varphi (-\ln  |x|^2  )   \\
&	- & 2 \eta  (|x|) \left ( x_1^2 + x_2^2 \right )   |x|^{2t-2} \varphi ' (-\ln  |x|^2  )    +  2t(2t-2) \eta  (|x|)   x_1 ^2 x_2^2 |x|^{2t-4} \varphi (-\ln  |x|^2  ) \nonumber \\
& - &  4t \eta  (|x|)   x_1^2 x_2^2 |x|^{2t-4} \varphi '  (-\ln |x|^2  )   -  4\eta '  (|x|)   x_1 ^2 x_2^2 |x|^{2t-3} \varphi ' (-\ln  |x|^2  )   \nonumber \\
			&  - & 2(2t-2)\eta   (|x|)   x_1 ^2 x_2^2 |x|^{2t-4} \varphi ' (-\ln  |x|^2  )  +4 \eta  (|x|)   x_1^2 x_2^2  |x|^{2t-4} \varphi '' (-\ln |x|^2  ) \nonumber
\end{eqnarray}

\vspace{.05in}
\noindent
For each $k$, choose $x^{(k)} \in \mathbb{R}^n$ such that
\begin{equation}
	\label{eqn:defn_x_(k)}
	\frac{x^{(k)} - R_k\zeta _0}{r_k} = \left (  e^{-\frac{1}{4t_k}}, 0, ... 0 \right ). \nonumber 
\end{equation}

\noindent
Then $x^{(k)}$ is in the ball $B_k$ and 

\begin{equation}
	\label{eqn:u_second_deriv_sample_points}
	\displaystyle \frac{\partial ^2 u}{\partial x_1 \partial x_2} \left (x^{(k)} \right ) \,\, = \,\,  \epsilon _k  \frac{\partial^2 u_{t_k} }{\partial x_1 \partial x_2}  \left ( \frac{x^{(k)}-R_k\zeta_0}{r_k} \right )  \,\, = \,\, \epsilon _k  \frac{\partial^2 u_{t_k} }{\partial x_1 \partial x_2}  \left (\left (  e^{-\frac{1}{4t_k}}, 0, ... 0 \right )\right ) .
\end{equation}

\noindent
Because $t_k < \frac{1}{4}$, we know $e^{-\frac{1}{4t_k}} < e^{-1} < \frac{1}{2}$, so in a neighborhood of $e^{-\frac{1}{4t_k}}$, $\eta(|x|) \equiv 1$ and $ \eta '(|x|)  =  \eta '' (|x|)  =0. $ Also note that $x_2=0$ at the point $\left (  e^{-\frac{1}{4t_k}}, 0, ... 0 \right )$, then by (\ref{eqn:2nd_derivatives_u_t_x_1_x_2}) we have
\allowdisplaybreaks
\begin{eqnarray}
& & 	\label{eqn: u_t_second_deriv_evaluate}
	\frac{\partial ^2 u_{t_k}}{\partial x_1 \partial x_2} \left (\left (  e^{-\frac{1}{4t_k}}, 0, ... 0 \right ) \right ) \nonumber \\
	 & = &  \left ( e^{-\frac{1}{4t_k}} \right )^{2t_k} \varphi  \left (-\ln \left ( e^{-\frac{1}{4t_k}} \right )^2 \right )   +  2t_k \left ( e^{-\frac{1}{4t_k}} \right )^{2t_k} \varphi  \left (-\ln \left ( e^{-\frac{1}{4t_k}} \right )^2 \right ) \nonumber \\
	  & - & 2 \left ( e^{-\frac{1}{4t_k}} \right )^{2t_k} \varphi ' \left (-\ln \left ( e^{-\frac{1}{4t_k}} \right )^2 \right ) \nonumber  \\
	& = & e^{-\frac{1}{2}}  \varphi \left ( \frac{1}{2t_k} \right ) +  2t_k \left ( e^{-\frac{1}{4t_k}} \right )^{2t_k} \varphi  \left (-\ln \left ( e^{-\frac{1}{4t_k}} \right )^2 \right ) - 2 e^{-\frac{1}{2}} \varphi ' \left ( \frac{1}{2t_k} \right ),   \nonumber 
\end{eqnarray} 

\vspace{.05in}
\noindent
where we purposefully did not simplify the second term.
Then (\ref{eqn:u_second_deriv_sample_points}) becomes
	\allowdisplaybreaks
\begin{eqnarray}
	\label{eqn: u_second_derivative_evaluate}
& &	 \frac{\partial ^2 u}{\partial x_1 \partial x_2} \left (x^{(k)} \right ) \\
	 & = & \epsilon_k e^{-\frac{1}{2}}  \varphi \left ( \frac{1}{2t_k} \right ) + \epsilon_k\cdot 2t_k \left ( e^{-\frac{1}{4t_k}} \right )^{2t_k} \varphi  \left (-\ln \left ( e^{-\frac{1}{4t_k}} \right )^2 \right ) - 2 \epsilon_k e^{-\frac{1}{2}} \varphi ' \left ( \frac{1}{2t_k} \right ) . \nonumber
\end{eqnarray}

\vspace{.05in}
\noindent
The second term in (\ref{eqn: u_second_derivative_evaluate}) goes to $0$ because 
$$  \epsilon_k \left \vert 2t_k \left ( e^{-\frac{1}{4t_k}} \right )^{2t_k} \varphi  \left (-\ln \left ( e^{-\frac{1}{4t_k}} \right )^2 \right ) \right \vert \leq \epsilon_k B_1  $$ by Lemma \ref{lemma:t_times_x^t_time_varphi}.
The third term in (\ref{eqn: u_second_derivative_evaluate}) goes to $0$ because  $ \frac{1}{2t_k} \to \infty$ and  $\displaystyle \lim_{s \to \infty} \varphi '(s)=0$.

\noindent
Now choose 
\begin{equation}
	\label{defn:epsilon_n}
	\epsilon _k = \frac{1}{\sqrt{\varphi \left ( \frac{1}{2t_k} \right ) }}. \nonumber
\end{equation}

\noindent
The first term in (\ref{eqn: u_second_derivative_evaluate}) becomes $$  \epsilon_k e^{-\frac{1}{2}}  \varphi \left ( \frac{1}{2t_k} \right ) =   e^{-\frac{1}{2}} \sqrt{ \varphi \left ( \frac{1}{2t_k} \right )} \to \infty \hspace{.2in} \text{as} \,\, k \to \infty. $$ Therefore,

  $$\displaystyle \lim _{k \to \infty} \frac{\partial ^2 u}{\partial x_1 \partial x_2} \left (x^{(k)} \right )  = \infty.$$  This shows that  $ \frac{\partial ^2 u}{\partial x_1 \partial x_2}$ is not bounded near the origin, and the lemma is proved.

\stop

\vspace{.2in}

\section{Construction for Theorem \ref{thm: kth_order}}
\label{section:higher_order}

\vspace{.1in}
The idea for constructing higher order examples for Theorem \ref{thm: kth_order} is the same as that for Theorem \ref{thm:main}, and we only need to replace $x_1x_2$ by the real or imaginary part of $(x_1+ix_2)^{k+2}$, where $k \in \mathbb{N}$.  For example, if $k=1$, then 
$$  (x_1+ix_2)^{3}  =\left (  x_1^3 - 3x_1x_2^2  \right ) + i \left ( 3x_1^2 x_2 -x_2^3  \right ), $$ so we may use either $x_1^3 - 3x_1x_2^2$ or $3x_1^2x_2-x_2^3$ in the construction.  For general $k$,
$$  (x_1+ix_2)^{k+2} = \sum _{l=0}^{k+2}  \binom{k+2}{l}  x_1^{k+2-l} \left ( ix_2 \right )  ^l. $$
 Evidently, the expressions for its real and imaginary parts are inconvenient to compute.  Thus to simplify the calculations we use complex variables for the first two components of $x$: for any $x=(x_1, x_2, x_3,...,x_n) \in \mathbb{R}^n$, denote
 \begin{equation*}
 	\label{eqn:define_z}
 	z=x_1+ix_2 \,\, \,\, \text{and} \,\, \,\, \bar{z}=x_1-ix_2.
 \end{equation*} 
Then
$$ x_1^2+x_2^2=z\bar{z}, \hspace{.4in} |x|^2=z\bar{z}+\sum_{j=3}^n x_j ^2, \hspace{.4in} \text{and} \hspace{.4in} \frac{\partial ^2 }{\partial x_1^2} + 	\frac{\partial ^2 }{\partial x_2^2} =4	\frac{\partial ^2 }{\partial \bar{z} \partial z},$$ and consequently
\begin{eqnarray*}
	\frac{\partial |x|}{ \partial z} = \frac{\bar{z}}{2|x|}, \hspace{.4in} \frac{\partial |x|}{ \partial \bar{z}} = \frac{z}{2|x|}, \hspace{.4in} \text{and} \hspace{.4in} \frac{\partial |x|}{ \partial x_j} = \frac{x_j}{|x|} \,\, \,\, (\text{when} \,\, j \geq 3).
\end{eqnarray*}  
 
 \vspace{.05in}
 \noindent
 Our strategy is to create a complex-valued function such that it is $(k+2)$-times differentiable in $\mathbb{R}^n$, its Laplacian is $C^k$, but $D^{k+2}u$ is unbounded. The real and imaginary parts of $u$ are two real-valued functions, and at least one of them would be a desired function that satisfies all the conditions in Theorem \ref{thm: kth_order}.  The proof is similar to that for Theorem \ref{thm:main}, so we will only present the key calculations and noticeable differences without repeating the entire proof.
 
\vspace{.05in}
 \noindent
 Recall that in the higher order case $\varphi(s)$ needs to be $(k+2)$-times differentiable and satisfy
 \begin{equation}
 	\label{eqn:higher_condition_phi}
 		\lim _{s \to \infty} \varphi  (s) = \infty,  \hspace{.2in} \lim _{s \to \infty} \varphi ' (s) = \cdots = \lim _{s \to \infty} \varphi ^{(k+2)}(s) = 0.
 \end{equation}
 
 
\vspace{.05in}
 \noindent
The building block function in this case needs to be modified into $$v(x) = z^{k+2}\varphi  (-\ln  |x|^2  ).$$ 

\noindent
The Laplacian of $v$ is
\begin{eqnarray*}
	\Delta v & = & 4	\frac{\partial ^2 v}{\partial \bar{z} \partial z} + \sum _{j=3}^n 	\frac{\partial ^2 v}{ \partial x_j^2}  \\
	& = & - (2n+4k)z^{k+2}|x|^{-2} \varphi '  (-\ln  |x|^2  )  + 4z^{k+2}|x|^{-2} \varphi ''  (-\ln  |x|^2  ).
\end{eqnarray*}

\vspace{.05in}
\noindent
It can be verified that $\Delta v$ is $C^k$, the partial derivative $	\dfrac{\partial ^{k+2} v}{ \partial z^{k+2}} $ is unbounded, and $v$ is not $(k+2)$-times differentiable.  Because this fact is not to be used in our constructions, we will not verify it here.


\vspace{.05in}
As in Section \ref{section:example_unbounded_hessian}, the next step is to  smooth out the function $v$.  Define $u_t: \mathbb{R}^n \to \mathbb{R} \,\, (n \geq 2)$ by
\begin{equation}
	\label{eqn:higher_order_defn_u_t}
	u_t(x) = 	\left\{
	\begin{array}{l l l}      
		0 &   \hspace{.1in} x=0, \\
		\noalign{\medskip}
		\eta(|x|)z^{k+2}|x|^{2t}\varphi  (-\ln  |x|^2  ) & \hspace{.1in} 0 < |x|< 1, \\
		\noalign{\medskip}
		0 &	\hspace{.1in} |x| \geq 1,
	\end{array}\right.
\end{equation}
where $\eta$ is the same as in (\ref{eqn: condition_eta}) and $\varphi$ satisfies (\ref{eqn:higher_condition_phi}).  This $u_t$ is a complex-valued $C^{k+2}$ function. The following is a key fact that will be used later.

\vspace{.05in}
\begin{lemma}
	\label{lemma: higher_u_t}
	For $u_t$ defined by (\ref{eqn:higher_order_defn_u_t}), the $k$-th partial derivatives of $\Delta u_t$ are all bounded by a constant independent of $t$.
\end{lemma}





\pf
The Laplacian of $u_t$ is  $$ \Delta u_t= 4\frac{\partial ^2 u_t}{\partial \bar{z} \partial z}  + \sum_{j=3}^{n} \frac{\partial ^2 u_t}{\partial x_j^2},    $$ where
\vspace{.05in}
\allowdisplaybreaks
\begin{eqnarray}
	\label{eqn:2nd_derivatives_u_t_z_z-bar}
	\frac{\partial ^2 u_t}{\partial \bar{z} \partial z} (x) 
	& = &  \frac{1}{4} \eta '' (|x|)  z^{k+3} \bar {z}|x|^{2t-2} \varphi  (-\ln  |x|^2   ) +   \frac{k+3}{2} \eta ' (|x|)  z^{k+2} |x|^{2t-1} \varphi  (-\ln  |x|^2   ) \nonumber \\
	& + & \frac{4t-1}{4}  \eta ' (|x|)   z^{k+3} \bar{z} |x|^{2t-3} \varphi  (-\ln  |x|^2   )
	-  \eta ' (|x|)   z^{k+3} \bar{z} |x|^{2t-3} \varphi  (-\ln  |x|^2  ) \nonumber \\
	&	+ & (k+3)t \eta  (|x|)   z^{k+2} |x|^{2t-2} \varphi  (-\ln  |x|^2  )  -  (k+2) \eta  (|x|)   z^{k+2} |x|^{2t-2} \varphi '  (-\ln |x|^2  )  \\ 
	& + & t(t-1) \eta  (|x|)   z^{k+3} \bar{z} |x|^{2t-4} \varphi  (-\ln |x|^2  )  -   (2t-1) \eta  (|x|)   z^{k+3} \bar{z} |x|^{2t-4} \varphi '  (-\ln ( |x|^2  ) \nonumber \\
	& - &   \eta '  (|x|)  z^{k+2} |x|^{2t-2} \varphi '  (-\ln  |x|^2  )  + \eta  (|x|)  z^{k+3} \bar{z} |x|^{2t-4} \varphi ''  (-\ln  |x|^2  ), \nonumber 
\end{eqnarray}
and
\allowdisplaybreaks
\begin{eqnarray}
	\label{eqn:2nd_derivatives_u_t_x_higher}
	\frac{\partial ^2 u_t}{\partial x_j^2} (x) & = &   \eta '' (|x|)   z^{k+2} x_j^2 |x|^{2t-2} \varphi  (-\ln  |x|^2   )   + \eta ' (|x|)  z^{k+2} |x|^{2t-1} \varphi  (-\ln  |x|^2   ) \nonumber  \\
	& + & (4t-1)  \eta ' (|x|)   z^{k+2} x_j^2|x|^{2t-3} \varphi  (-\ln  |x|^2   )
	-4  \eta ' (|x|)  z^{k+2} x_j^2 |x|^{2t-3} \varphi '  (-\ln  |x|^2  )  \nonumber \\
	& +  & 	 2t \eta  (|x|)   z^{k+2} |x|^{2t-2} \varphi (-\ln  |x|^2 )  +  2t(2t-2) \eta  (|x|)  z^{k+2} x_j^2 |x|^{2t-4} \varphi (-\ln  |x|^2  ) \\ 
	& + & (4-8t) \eta  (|x|)   z^{k+2} x_j^2 |x|^{2t-4} \varphi ' (-\ln  |x|^2  )  -     2 \eta   (|x|)   z^{k+2}  |x|^{2t-2} \varphi '  (-\ln  |x|^2 ) \nonumber \\
	&  + &   4 \eta  (|x|)   z^{k+2} x_j^2 |x|^{2t-4} \varphi ''  (-\ln |x|^2  ). \nonumber
\end{eqnarray}

\vspace{.05in}
\noindent
We will show that all the $k$-th partial derivatives of each term in (\ref{eqn:2nd_derivatives_u_t_z_z-bar}) and (\ref{eqn:2nd_derivatives_u_t_x_higher}) are bounded by a constant independent of $t$.  In the subsequent discussions in this section, we will use $C$ to denote a constant that depends on $\eta$, $\varphi$, $n$, $k$ and is independent of $t$.


\vspace{.05in}
\noindent
We start with (\ref{eqn:2nd_derivatives_u_t_z_z-bar}). The first term of (\ref{eqn:2nd_derivatives_u_t_z_z-bar}) is bounded by 
\begin{equation}
	\label{eqn:partial_z_z_bar_first_term}
	 \frac{1}{4} \eta '' (|x|)  z^{k+3} \bar {z}|x|^{2t-2} \left | \varphi  (-\ln  |x|^2   ) \right | \leq C|x|^{k+2} \left | \varphi  (-\ln  |x|^2   ) \right |.
\end{equation}

\vspace{.05in}
\noindent
Its first derivatives may be taken with respect to $z$, $\bar{z}$, or $x_j \,\, (j \geq 3)$.

\vspace{.05in}
\begin{itemize}

	\item If we take its derivative with respect to $z$, then
	\begin{eqnarray}
	& &	\frac{\partial}{\partial z} \left (  \eta '' (|x|)  z^{k+3} \bar {z}|x|^{2t-2} \varphi  (-\ln  |x|^2   )     \right ) \nonumber \\
		& = & \left ( \eta ^{(3)} (|x|) \frac{\bar{z}}{2|x|} \right ) z^{k+3} \bar {z}|x|^{2t-2} \varphi  (-\ln  |x|^2   )    +  \eta '' (|x|) \left ( (k+3)z^{k+2} \right ) \bar {z}|x|^{2t-2} \varphi  (-\ln  |x|^2   )  \nonumber \\
		& + & \eta '' (|x|)  z^{k+3} \bar {z} \left ( (2t-2)|x|^{2t-3}\frac{\bar{z}}{2|x|} \right ) \varphi  (-\ln  |x|^2   ) +  \eta '' (|x|)  z^{k+3} \bar {z}|x|^{2t-2} \left ( \varphi ' (-\ln  |x|^2   ) \frac{-\bar{z}}{|x|^2} \right ). \nonumber
	\end{eqnarray}

\vspace{.05in}
\noindent
Since $|z| \leq |x| \leq 1$ and $|\bar{z}| \leq |x| \leq 1$, the first term is bounded by 
\begin{equation}
	\label{eqn:z_derivative_1st_term_in_partial_z_z_bar}
	C|x|^{k+2+2t} \left | \varphi   (-\ln  |x|^2 ) \right |,
\end{equation} 
the second and third terms are bounded by 
\begin{equation}
	\label{eqn:z_derivative_2nd_3rd_term_in_partial_z_z_bar}
	C|x|^{k+1+2t} \left | \varphi   (-\ln  |x|^2 ) \right |,
\end{equation}
 the fourth term is bounded by 
 \begin{equation}
 	\label{eqn:z_derivative_4th_term_in_partial_z_z_bar}
 	C|x|^{k+1+2t} \left | \varphi '  (-\ln  |x|^2 ) \right |.
 \end{equation}

\vspace{.05in}
\noindent
 Note that both (\ref{eqn:z_derivative_1st_term_in_partial_z_z_bar}) and (\ref{eqn:z_derivative_2nd_3rd_term_in_partial_z_z_bar}) are bounded by $C|x|^{k+1} \left | \varphi   (-\ln  |x|^2 )\right |, $ and (\ref{eqn:z_derivative_4th_term_in_partial_z_z_bar}) is bounded by $C|x|^{k+1} \left |\varphi '  (-\ln  |x|^2 )\right |. $  
 Therefore, this derivative is bounded by 
 $$ C|x|^{k+1} \left | \varphi   (-\ln  |x|^2 ) | + C|x|^{k+1} | \varphi '  (-\ln  |x|^2 ) \right | . $$
 
\vspace{.05in}
	\item If we take its derivative with respect to $\bar{z}$, then
\begin{eqnarray}
	& &	\frac{\partial}{\partial \bar{z}} \left (  \eta '' (|x|)  z^{k+3} \bar {z}|x|^{2t-2} \varphi  (-\ln  |x|^2   )     \right ) \nonumber \\
	& = & \left ( \eta ^{(3)} (|x|) \frac{z}{2|x|} \right ) z^{k+3} \bar {z}|x|^{2t-2} \varphi  (-\ln  |x|^2   )    +  \eta '' (|x|) z^{k+3} |x|^{2t-2} \varphi  (-\ln  |x|^2   )  \nonumber \\
	& + & \eta '' (|x|)  z^{k+3} \bar {z} \left ( (2t-2)|x|^{2t-3}\frac{z}{2|x|} \right ) \varphi  (-\ln  |x|^2   ) +  \eta '' (|x|)  z^{k+3} \bar {z}|x|^{2t-2} \left ( \varphi ' (-\ln  |x|^2   ) \frac{-z}{|x|^2} \right ). \nonumber
\end{eqnarray}

\vspace{.05in}
\noindent
By similar argument we know that this derivative is also bounded by 
 $$ C|x|^{k+1} \left | \varphi   (-\ln  |x|^2 ) \right | + C|x|^{k+1} \left | \varphi '  (-\ln  |x|^2 ) \right |. $$
 
\vspace{.05in}
	\item If we take its derivative with respect to $x_j$, then
\begin{eqnarray}
	& &	\frac{\partial}{\partial x_j} \left (  \eta '' (|x|)  z^{k+3} \bar {z}|x|^{2t-2} \varphi  (-\ln  |x|^2   )     \right ) \nonumber \\
	& = & \left ( \eta ^{(3)} (|x|) \frac{x_j}{|x|} \right ) z^{k+3} \bar {z}|x|^{2t-2} \varphi  (-\ln  |x|^2   )   +  \eta '' (|x|)  z^{k+3} \bar {z} \left ( (2t-2)|x|^{2t-3}\frac{x_j}{|x|} \right ) \varphi  (-\ln  |x|^2   ) \nonumber \\
	& + &  \eta '' (|x|)  z^{k+3} \bar {z}|x|^{2t-2} \left ( \varphi ' (-\ln  |x|^2   ) \frac{-2x_j}{|x|^2} \right ). \nonumber
\end{eqnarray}

\vspace{.05in}
\noindent
Again, this derivative is bounded by 
$$ C|x|^{k+1} \left | \varphi   (-\ln  |x|^2 ) \right | + C|x|^{k+1} \left |  \varphi '  (-\ln  |x|^2 ) \right |. $$
\end{itemize}

\vspace{.05in}
\noindent
In conclusion, regardless of which variable we differentiate with, the first partial derivative of $ \eta '' (|x|)  z^{k+3} \bar {z}|x|^{2t-2} \varphi  (-\ln  |x|^2   )      $ is bounded by 
$$ C|x|^{k+1} \left |\varphi   (-\ln  |x|^2 ) \right | + C|x|^{k+1} \left | \varphi '  (-\ln  |x|^2 ) \right |. $$
Comparing the first term of this bound, $C|x|^{k+1} \left | \varphi   (-\ln  |x|^2 ) \right |$, to the right hand side of (\ref{eqn:partial_z_z_bar_first_term}), we see that the power of $|x|$ decreased from $k+2$ to $k+1$.  By the same type of  calculations, the $k$-th partial derivative of  $ \eta '' (|x|)  z^{k+3} \bar {z}|x|^{2t-2} \varphi  (-\ln  |x|^2) $ will be bounded by  
$$ C|x|^{2} \left |\varphi   (-\ln  |x|^2 ) \right | + \sum _{l=1}^k C|x|^{k+2-l} \left |\varphi ^{(l)}  (-\ln  |x|^2 ) \right |. $$

\noindent
Because $$ \lim _{s \to \infty} \varphi ' (s) = \cdots = \lim _{s \to \infty} \varphi ^{(k+2)}(s) = 0, $$ $\displaystyle \sum _{l=1}^k C|x|^{k+2-l} \left | \varphi ^{(l)}  (-\ln  |x|^2 ) \right |$ is bounded by a constant independent of $t$.  By Lemma \ref{lemma:x-times-varphi}, $\displaystyle  C|x|^{2} \left | \varphi   (-\ln  |x|^2 ) \right |$ is also bounded by a constant independent of $t$.  Therefore, the $k$-th derivative of the first term of (\ref{eqn:2nd_derivatives_u_t_z_z-bar}) is bounded by a constant independent of $t$.  

\vspace{.05in}
\noindent
All the other terms in (\ref{eqn:2nd_derivatives_u_t_z_z-bar}) and (\ref{eqn:2nd_derivatives_u_t_x_higher}) can be estimated in the same way.  Only the 5th term in (\ref{eqn:2nd_derivatives_u_t_z_z-bar}), $(k+3)t \eta  (|x|)   z^{k+2} |x|^{2t-2} \varphi  (-\ln  |x|^2  ) $, and the 5th term in (\ref{eqn:2nd_derivatives_u_t_x_higher}), $  2t \eta  (|x|)   z^{k+2} |x|^{2t-2} \varphi (-\ln  |x|^2 ) $, are slightly different. By the same calculations as above, the $k$-th derivatives of  $$(k+3)t \eta  (|x|)   z^{k+2} |x|^{2t-2} \varphi  (-\ln  |x|^2  ) $$ are bounded by 
$$ Ct|x|^{2t} \left | \varphi   (-\ln  |x|^2 ) \right | + \sum _{l=1}^k Ct|x|^{k-l} \left | \varphi ^{(l)}  (-\ln  |x|^2 ) \right | . $$
The $\displaystyle \sum _{l=1}^kCt|x|^{k-l} \left | \varphi ^{(l)}  (-\ln  |x|^2 ) \right | $ term is bounded by a constant independent of $t$, by the same argument as above.  To estimate  $ Ct|x|^{2t} \left | \varphi   (-\ln  |x|^2 ) \right | $ we need to invoke Lemma \ref{lemma:t_times_x^t_time_varphi} instead of Lemma \ref{lemma:x-times-varphi}: since
$$ 2t|x|^{2t} \left | \varphi   (-\ln  |x|^2 ) \right | \leq B_1, $$
 $ Ct|x|^{2t} \left | \varphi   (-\ln  |x|^2 ) \right | $ is bounded by a constant independent of $t$.  Consequently, the $k$-th derivatives of  $$(k+3)t \eta  (|x|)   z^{k+2} |x|^{2t-2} \varphi  (-\ln  |x|^2  ) $$ are bounded by a constant independent of $t$.  The $k$-th derivatives of the 5th term in (\ref{eqn:2nd_derivatives_u_t_x_higher}), $  2t \eta  (|x|)   z^{k+2} |x|^{2t-2} \varphi (-\ln  |x|^2 ) $, can be handled in the same way. This completes the proof of the lemma.
 
 \stop

\vspace{.1in}
Then we define $u$ by $$ u(x) = \sum _{l=1}^{\infty} \epsilon_l r_l ^{k+2} u_{t_l} \left ( \frac{x-R_l \zeta_0}{r_l} \right ) ,$$ where $R_l$, $r_l$, $\zeta_0$, $t_l$, and $\epsilon_l$ are the same as in Section \ref{section:example_unbounded_hessian}; namely, $R_l$ and $r_l$ are decreasing sequences and  $$R_l -r_l > R_{l+1}+r_{l+1}, \,\, \,\,\,\, \zeta_0 =\left ( \frac{1}{\sqrt{2}}, \frac{1}{\sqrt{2}}, ..., \frac{1}{\sqrt{2}}  \right ), \,\,\,\,\,\, \displaystyle \lim_{l \to \infty} t_l=0, \,\,\,\,  \,\,	\epsilon _l = \frac{1}{\sqrt{\varphi \left ( \frac{1}{2t_l} \right ) }}.$$ Thus by the same argument as in Section \ref{section:example_unbounded_hessian} we know that this infinite sum actually only has a single term for any given $x$ value. 

\vspace{.05in}
We first show that $\Delta u$ is $C^k$. Note that here the power of $r_l$ is $k+2$ as opposed to 2 in Section \ref{section:example_unbounded_hessian}, so $$ \Delta u(x) = \sum _{l=1}^{\infty} \epsilon_l r_l^{k}\Delta u_{t_l} \left ( \frac{x-R_l \zeta_0}{r_l} \right ).  $$ 
By Lemma \ref{lemma: higher_u_t}, the $k$-th derivatives of $\Delta u_{t_l} $ are uniformly bounded by a constant independent of $t_l$.  Then since $ \displaystyle \lim_{l \to \infty} \epsilon_l = 0$, we conclude that the $k$-th derivatives of $\Delta u$ all approach $0$ as $|x| \to 0$.   Recall that by construction  $u = 0$ on all of the coordinate hyperplanes, so all partial derivatives of $u$ of any order is $0$ at the origin.  In particular, all the $k$-th derivatives of $\Delta u$ at the origin is 0. Therefore, all the $k$-th derivatives of $\Delta u$ are continuous at the origin, and consequently $\Delta u$ is $C^k$ throughout $\mathbb{R}^n$.

\vspace{.05in}
Next, we show that some of the $(k+2)$-th derivatives of $u$ is unbounded.  Precisely, we will show that $\dfrac{\partial ^{k+2} u}{\partial z^{k+2}}$ is unbounded.
We start with a close look at the first and second partial derivatives of $u_t$ with respect to $z$.  

\vspace{.05in}
\noindent
Note that the power of $z$ in $$ u_t = 	\eta(|x|)z^{k+2}|x|^{2t}\varphi  (-\ln  |x|^2  ) $$ is $k+2$.  After one differentiation with respect to $z$, one of the terms in its derivative is $$ 	(k+1)\eta(|x|)z^{k+1}|x|^{2t}\varphi  (-\ln  |x|^2  ), $$ which is the first term in the following formula for $	\frac{\partial  u_t}{ \partial z}$:

\allowdisplaybreaks
\begin{eqnarray}
	\label{eqn:higher_u_t_derivative_z}
		\frac{\partial u_t}{\partial z} (x) & = & (k+2) \eta (|x|) z^{k+1}|x|^{2t}  \varphi  (-\ln  |x|^2  ) + \frac{1}{2} \eta ' (|x|) z^{k+2}\bar{z}|x|^{2t-1} \varphi  (-\ln  |x|^2  ) \nonumber \\
		&  + & t \eta (|x|) z^{k+2}\bar{z} |x|^{2t-2} \varphi  (-\ln  |x|^2  )   - \eta (|x|)z^{k+2}\bar{z}|x|^{2t-2} \varphi '  (-\ln  |x|^2  ).  
		\end{eqnarray}

\vspace{.05in}
\noindent
After another differentiation with respect to $z$, there will be one term, $$(k+1)(k+2) \eta  (|x|)   z^{k} |x|^{2t} \varphi  (-\ln  |x|^2  ).  $$ where the power of $z$ is $k$.  That is the first term in the following formula for $	\frac{\partial ^2 u_t}{ \partial z^2}$:
\vspace{.05in}
\allowdisplaybreaks	
	\begin{eqnarray}
			\label{eqn:higher_u_t_derivative_z^2}
				\frac{\partial ^2 u_t}{ \partial z^2} (x) 
			& = & (k+1)(k+2) \eta  (|x|)   z^{k} |x|^{2t} \varphi  (-\ln  |x|^2  ) + \frac{1}{4} \eta '' (|x|)  z^{k+2} \bar {z}^2|x|^{2t-2} \varphi  (-\ln  |x|^2   ) \nonumber \\
			& + &  (k+2) \eta ' (|x|)  z^{k+1} \bar{z} |x|^{2t-1} \varphi  (-\ln  |x|^2   )  +  \frac{4t-1}{2}  \eta ' (|x|)   z^{k+2} \bar{z}^2 |x|^{2t-3} \varphi  (-\ln  |x|^2   ) \nonumber \\
		&	+ &  2t  (k+2) \eta  (|x|)   z^{k+1} \bar{z} |x|^{2t-2} \varphi   (-\ln |x|^2  ) + t  (t-1) \eta  (|x|)   z^{k+2} \bar{z}^2 |x|^{2t-4} \varphi  (-\ln ( |x|^2  )  \\
			 &	- & \eta ' (|x|)   z^{k+2} \bar{z}^2 |x|^{2t-3} \varphi ' (-\ln  |x|^2  ) -  2(k+2) \eta  (|x|)   z^{k+1} \bar{z} |x|^{2t-2} \varphi ' (-\ln |x|^2  )  \nonumber \\
			& - &  (2t-1) \eta   (|x|)  z^{k+2}\bar{z}^2 |x|^{2t-4} \varphi '  (-\ln  |x|^2  )  + \eta  (|x|)  z^{k+2} \bar{z}^2 |x|^{2t-4} \varphi ''  (-\ln  |x|^2  ), \nonumber 
\end{eqnarray}

\vspace{.05in}
\noindent
After $(k+2)$-times of differentiation with respect to $z$,  one of the terms in $	\dfrac{\partial ^{k+2} u_t}{ \partial z^{k+2}}  $ is  $$(k+2)! \eta  (|x|)  |x|^{2t} \varphi  (-\ln  |x|^2  ). $$  As will be shown later, this term is crucial to proving that $\dfrac{\partial ^{k+2} u}{\partial z^{k+2}}$ is unbounded.

\vspace{.05in}
\noindent
For each $l$, as we did in Section \ref{section:example_unbounded_hessian}, choose $x^{(l)} \in \mathbb{R}^n$ such that 
$$
	\frac{x^{(l)} - R_l\zeta _0}{r_l} = \left (  e^{-\frac{1}{4t_l}}, 0, ... 0 \right ). 
$$
\noindent
As discussed in Section \ref{section:example_unbounded_hessian}, in a neighborhood of 
$\left (  e^{-\frac{1}{4t_l}}, 0, ... ,0 \right ),  $ $\eta (|x|) \equiv 1$ and $\eta '(|x|) = \eta ''(|x|)  = 0$.  Therefore, when we evaluate $	\frac{\partial  u_t}{ \partial z}$ and 	$\frac{\partial ^2 u_t}{ \partial z^2}$ in a neighborhood of $\left (  e^{-\frac{1}{4t_l}}, 0, ..., 0 \right ), $ all the terms in (\ref{eqn:higher_u_t_derivative_z}) and (\ref{eqn:higher_u_t_derivative_z^2}) that have an $\eta ' $ or $\eta '' $ factor will disappear.  For that reason in the discussion that follows, we will only consider the terms that have an $\eta$ factor.

\vspace{.05in}
\noindent
Then (\ref{eqn:higher_u_t_derivative_z}) becomes 
$$  (k+2)  z^{k+1}|x|^{2t_l}  \varphi  (-\ln  |x|^2  )  +  t_l  z^{k+2}\bar{z} |x|^{2t_l-2} \varphi  (-\ln  |x|^2  )   - z^{k+2}\bar{z}|x|^{2t_l-2} \varphi '  (-\ln  |x|^2  ).   $$

\noindent
Besides the first term, the remaining terms are bounded by 
$$ |x|^{k+1+2t_l} \varphi  (-\ln  |x|^2   ) + |x|^{k+1+2t_l} \varphi ' (-\ln  |x|^2   ).   $$

\noindent
And (\ref{eqn:higher_u_t_derivative_z^2}) becomes
\begin{eqnarray*}
& &	(k+1)(k+2)  z^{k} |x|^{2t_l} \varphi  (-\ln  |x|^2  ) +2t_l  (k+2)    z^{k+1} \bar{z} |x|^{2t_l-2} \varphi   (-\ln |x|^2  ) \\
	& + & t_l  (t_l-1)    z^{k+2} \bar{z}^2 |x|^{2t_l-4} \varphi  (-\ln ( |x|^2  )   -  2(k+2)    z^{k+1} \bar{z} |x|^{2t_l-2} \varphi ' (-\ln |x|^2  ) \\
	& - &  (2t_l-1)   z^{k+2}\bar{z}^2 |x|^{2t_l-4} \varphi '  (-\ln  |x|^2  )  +   z^{k+2} \bar{z}^2 |x|^{2t_l -4} \varphi ''  (-\ln  |x|^2  ) 
\end{eqnarray*}

\noindent
Besides the first term, the remaining terms are bounded by 
$$ Ct_l|x|^{k+2t_l} \varphi  (-\ln  |x|^2   ) + C|x|^{k+2t_l} \varphi ' (-\ln  |x|^2   ) + C|x|^{k+2t_l} \varphi '' (-\ln  |x|^2   ).   $$

\vspace{.05in}
\noindent
By the same process, in a neighborhood of 
$\left (  e^{-\frac{1}{4t_l}}, 0, ... 0 \right ),$  $\dfrac{\partial ^{k+2} u_{t_l}}{\partial z^{k+2}}$ is equal to $$ (k+2)!   |x|^{2t_l} \varphi  (-\ln  |x|^2  )   $$ plus some other terms that are bounded by
\begin{equation}
	\label{eqn:higher_example_bounded_terms}
  Ct_l|x|^{2t_l} \varphi  (-\ln  |x|^2  )  +  Ct_l|x|^{2t_l} \varphi ' (-\ln  |x|^2  ) + \cdots Ct_l|x|^{2t_l} \varphi ^{(k+2)} (-\ln  |x|^2  ).
\end{equation}

\vspace{.05in}
\noindent
By Lemma \ref{lemma:t_times_x^t_time_varphi} and the fact that $$ \lim _{s \to \infty} \varphi ' (s) = \cdots = \lim _{s \to \infty} \varphi ^{(k+2)}(s) = 0, $$ we know
	(\ref{eqn:higher_example_bounded_terms})
	  is bounded by a constant independent of $t_l$.

\vspace{.05in}
\noindent 
Now we look at 
\begin{eqnarray*}
	\frac{\partial ^{k+2} u}{ \partial z^{k+2}} (x^{(l)}) & = & \epsilon_l  \frac{\partial ^{k+2} u_{t_l}}{ \partial z^{k+2}} \big ( (e^{-\frac{1}{4t_l}}, 0, ..., 0 ) \big ).
\end{eqnarray*}
After evaluating (\ref{eqn:higher_example_bounded_terms}) at the point $\left ( e^{-\frac{1}{4t_l}}, 0, ..., 0 \right )$ and multiplying the result with $\epsilon_l  $, it goes to 0 as $\epsilon_l \to 0$. 

\noindent
The first term of $ 	\dfrac{\partial ^{k+2} u}{ \partial z^{k+2}} (x^{(l)})$ is equal to 
$$ \epsilon_l (k+2)!   \left ( e^{-\frac{1}{4t_l}} \right )^{2t_l} \varphi  \left (-\ln  \left ( e^{-\frac{1}{4t_l}} \right )^2  \right ) = \frac{(k+2)!}{\sqrt{e}} \epsilon_l \varphi \left ( \frac{1}{2t_l}\right ). $$
Recall that  $$\epsilon_l  = \frac{1}{\sqrt{\varphi \left ( \frac{1}{2t_l}\right ) }},  $$ hence $\epsilon_l \varphi \left ( \frac{1}{2t_l}\right ) \to \infty $ as $l \to \infty$. Consequently, $ 	\dfrac{\partial ^{k+2} u}{ \partial z^{k+2}} (x^{(l)}) \to \infty$ as $l \to \infty$, which implies that  $ 	\dfrac{\partial ^{k+2} u}{ \partial z^{k+2}} $ is unbounded near the origin.  Since $ \dfrac{\partial }{\partial z} = \dfrac{1}{2} \left (  \dfrac{\partial }{\partial x_1} - i  \dfrac{\partial }{\partial x_2} \right ), $ as a result we know that some of the partial derivatives of $u$ with respect to the $x_1$ and $x_2$ variables are unbounded near the origin.

\vspace{.05in}
Finally, we need to show that $ u $ is $(k+2)$-times differentiable at the origin.  Recall that because $u=0$ on all the coordinate hyperplanes, all partial derivatives of $u$ at the origin is 0, which include all the $(k+1)$-th partial derivatives of $u$.  Denote an arbitrary $(k+1)$-th partial derivative of $u$ as $D^{\gamma}u$, where $\gamma = (\gamma _1, ..., \gamma_ n)$ and $\gamma_1 + \cdots + \gamma_n=k+1$, then
\begin{eqnarray*}
 \frac{ \left | D^{\gamma}u(x) -D^{\gamma}u(0) \right | }{|x|} = \frac{ \left | D^{\gamma}u(x)  \right | }{|x|} & \leq & \frac{\epsilon_l r_l \left |  D^{\gamma}u_{t_l} \left ( \frac{x-R_l \zeta_0}{r_l}\right )  \right | }{R_l - \frac{2}{3}r_l}  \\
 & = & \frac{ \epsilon_l  \left |   D^{\gamma}u_{t_l} \left ( \frac{x-R_l \zeta_0}{r_l}\right )  \right | }{\frac{R_l}{r_l} - \frac{2}{3}} \\
 & < & 3\epsilon_l \left  |   D^{\gamma}u_{t_l} \left ( \frac{x-R_l \zeta_0}{r_l}\right )  \right |	,
\end{eqnarray*}
  where we have used the fact that $|x| \geq R_l - \frac{2}{3}r_l$ as proved in Section \ref{section:example_unbounded_hessian} and that $R_l > r_l$.
  
 \vspace{.05in}
  \noindent
  By calculations similar to those in this section, we can show that $\left |  D^{\gamma}u_{t_l} \left ( \frac{x-R_l \zeta_0}{r_l}\right )  \right |	$ is bounded by a constant independent of $t$.  As a result, $$ \frac{ \left | D^{\gamma}u(x) -D^{\gamma}u(0) \right | }{|x|} \to 0  \hspace{.4in} \text{as} \,\, \,\, |x| \to 0,  $$ which implies $u$ is $(k+2)$-times differentiable at the origin.

\vspace{.05in}
 Thus we can conclude that as a complex-valued function, $u$ is $(k+2)$-times differentiable at 0, $\Delta u$ is $C^k$ throughout $\mathbb{R}^n$, but $D^{k+2}u$ is unbounded near 0.  The real and imaginary parts of $u$ are two real-valued functions that are $(k+2)$-times differentiable at 0, their Laplacian are $C^k$ throughout $\mathbb{R}^n$, and at least one of them has some unbounded $(k+2)$-th partial derivatives.  Therefore, we have found a function that satisfies all the conditions in Theorem \ref{thm: kth_order}.
 
\vspace{.4in}

\section {Appendix: Proof of Proposition \ref{thm:gradient}}
\label{section:gradient_proof}

\vspace{.1in}

 The proof of Proposition \ref{thm:gradient} is based on a method that was introduced in \cite{X-Wang} and elaborated in detail in \cite{Oton}.  Note that after a translation we can assume $x$ or $y$ is at the origin, so we only need to prove that for $|z| < \frac{1}{16}$, (here the $z$ is a point in $\mathbb{R}^n$, not a complex variable as was used in the previous section), we have
\begin{equation}
	\label{eqn: first_derivative}
	\big | Du(z) - Du (0) \big |\leq C|z| \left ( \sup_{B_1} |u| + \sup_{B_1} |f| + \int _{|z|}^1 \frac{w(r)}{r} dr \right ) .  
\end{equation}
For $|z| \geq \frac{1}{16}$ the estimate is also true by a covering argument (see \cite{Oton}).  

\vspace{.05in}
\noindent
First, we recall three elementary estimates (see \cite{Oton}) that will be used frequently in this proof.

\vspace{.05in}
\noindent
If a function $v$ satisfies $\Delta v=0$ in $B_r$, then for any positive integer $k$,
\begin{equation}
	\label{eqn:estimate_Laplace}
	\| D^k v \|_{L^{\infty}(B_{\frac{r}{2}})} \leq Cr^{-k} \| v \|_{L^{\infty}(B_r)},
\end{equation}
where $C$ only depends on $n$ and $k$. 

\vspace{.05in}
\noindent
If a function $v$ satisfies $\Delta v=\lambda$ in $B_r$, where $\lambda$ is a constant and $r <1$, then 
\begin{equation}
	\label{eqn:estimate_Delta_with_constant}
	\| D v \|_{L^{\infty}(B_{\frac{r}{2}})} \leq C \left (  r^{-1} \| v \|_{L^{\infty}(B_r)} + r | \lambda | \right ).
\end{equation}

\vspace{.05in}
\noindent
If a function $v$ satisfies $\Delta v=f$ in $B_r$, where $f$ is a given bounded function, then the scaled maximum principle states that 
\begin{equation}
	\label{eqn:scaled_max}
	\| u \|_{L^{\infty}(B_r)} \leq \|u\|_{L^{\infty}(\partial B_r)} + Cr^2 \| f \|_{L^{\infty}(B_r)}.
\end{equation}

\vspace{.05in}
\noindent
Now we are ready to prove (\ref{eqn: first_derivative}).  For $k=0,1,2,...$, let $u_k$ be the solution to 

$$
\left\{
\begin{array}{r l l}      
	\displaystyle\Delta u_k & = & f(0)   \hspace{.2in} \text {in} \hspace{.1in} B_{2^{-k}} , \\
	\noalign{\smallskip}
	u_k & = & u   \hspace{.4in} \text {on} \hspace{.1in} \partial B_{2^{-k}}. 	
\end{array}\right.
$$

\vspace{.05in}
\noindent
Then $\displaystyle \Delta (u_k-u) = f(0)-f$ in $ B_{2^{-k}}$ and $u_k-u=0$ on $\partial B_{2^{-k}}$. By the scaled maximum principle it follows that
\begin{eqnarray}
	\label{eqn:u-k_minus_u}
	\| u_k-u \|_{L^{\infty}(B_{2^{-k}})} & \leq & C\left (2^{-2k} \right ) \|  f(0) -f \|_{L^{\infty}(B_{2^{-k}})}  \\
	& \leq & C \left ( 2^{-2k} \right ) \omega (2^{-k}), \nonumber
\end{eqnarray}

\noindent
and therefore
\begin{eqnarray}
	\label{eqn:u_K+1_minus_u_k}
	\| u_{k+1}-u_k \|_{L^{\infty}(B_{2^{-k-1}})} & \leq & \| u_{k+1}-u \|_{L^{\infty}(B_{2^{-k-1}})} + \| u_k-u \|_{L^{\infty}(B_{2^{-k}})}  \\
	& \leq &   C \left (2^{-2(k+1)} \right ) \omega  ( 2^{-(k+1)}  ) + C \left (2^{-2k} \right ) \omega (2^{-k}) \nonumber \\
	& \leq & C \left (2^{-2k} \right ) \omega (2^{-k}) . \nonumber
\end{eqnarray}

\noindent
Then since $u_{k+1} - u_k$ is harmonic, by (\ref{eqn:estimate_Laplace}) we have
\begin{eqnarray}
	\label{eqn:u_y+1-u_k}
	\|  Du_{k+1} -Du_k \|_{L^{\infty}(B_{2^{-k-2}})} & \leq & C \left (2^{k+1} \right ) \| u_{k+1}-u_k \|_{L^{\infty}(B_{2^{-k-1}})} \\
	& \leq & C \left ( 2^{-k} \right ) \omega (2^{-k}).  \nonumber
\end{eqnarray}

\noindent
For any $|z| \leq \frac{1}{16}$, choose $k \in \mathbb{N}$ such that $$ 2^{-k-4} \leq |z| \leq 2^{-k-3}. $$  We will estimate $\left \vert Du(z) - Du(0)  \right \vert$ by
\begin{equation}
	\label{eqn: 3_terms}
	\left \vert   Du(z) - Du(0)  \right \vert \leq
	\left \vert   Du(0) - Du_k(0)  \right \vert
	+
	\left \vert   Du(z) - Du_k(z)  \right \vert  + 	\left \vert   Du_k(z) - Du_k(0)  \right \vert.
\end{equation}

\vspace{.05in}
\noindent
We are going to estimate these three terms separately. First, we claim that 
\begin{equation*}
	\lim _{k \to \infty} Du_k(0)=Du(0).
\end{equation*}

\vspace{.05in}
\noindent
To see this, let $\tilde{u}(x)= u(0)+x\cdot Du(0)$ be the linear approximation of $u$ at 0.  Then $Du(0)=D\tilde{u}(0)$ and $|\tilde{u}(x)-u(x)|=o(|x|)$.  Thus
\begin{eqnarray*}
	& &\left \vert Du_k(0) - Du(0)   \right \vert \\
	 &  \leq &  \|  Du_k - D\tilde{u} \|_{L^{\infty}(B_{2^{-k-1}})} \\
	& \leq & C \left ( 2^k \right )  \|  u_k - \tilde{u} \|_{L^{\infty}(B_{2^{-k}})} \hspace{1.2in} \big (\text{ by (\ref{eqn:estimate_Laplace}) \big )}\\
	& = & C \left ( 2^k \right )  \left ( \|  u_k - \tilde{u} \|_{L^{\infty}(\partial B_{2^{-k}})} + 2^{-2k}|f(0)| \right )  \hspace{.1in} \left ( \text{apply} \hspace{.05in} (\ref{eqn:scaled_max}) \,\, \text{to} \,\, \Delta (u_k-\tilde{u}) = f(0) \right )\\
	& = & C \left ( 2^k \right )  \|  u - \tilde{u} \|_{L^{\infty}(\partial B_{2^{-k}})}  + C \left  ( 2^{-k} \right )|f(0)|\\
	& \leq & C \left ( 2^k \right ) \cdot o(2^{-k}) + C \left  ( 2^{-k} \right )|f(0)| \\
	& \to & 0  \hspace{.2in} \text{as}\hspace{.2in}  k \to \infty.
\end{eqnarray*}

\vspace{.05in}
\noindent
Then we can write $ \displaystyle Du_k(0) -Du(0)=\sum _{j=k}^{\infty}  \left ( Du_j(0)-Du_{j+1}(0) \right ), $ and consequently
\begin{eqnarray}
	\label{eqn:estimate_u_k-u_at_0}
	|Du_k(0) -Du(0)| & \leq & \sum _{j=k}^{\infty} |Du_j(0)-Du_{j+1}(0)|  \\
	& \leq & C \sum _{j=k}^{\infty} \left (2^{-j}\right ) \omega ( 2^{-j} ) \hspace{.6in} \big ( \text{by (\ref{eqn:u_y+1-u_k})} \big ) \nonumber \\
	& \leq & C \sum _{j=k}^{\infty} \left (2^{-j}\right ) \omega ( 2^{-k} ) \nonumber \\
	& = & C(2^{-k}) \omega ( 2^{-k} )  \nonumber \\
	& \leq & C|z| \omega ( 2^{-k} ). \nonumber
\end{eqnarray}

\vspace{.05in}
\noindent
Next, we estimate the term $ 	\left \vert   Du(z) - Du_k(z)  \right \vert$.  Let $v_j$ be the solution of 

$$
\left\{
\begin{array}{l l l}      
	\displaystyle\Delta v_j & = & f(z)   \hspace{.2in} \text {in} \hspace{.1in} B_{2^{-j}}(z) , \\
	\noalign{\smallskip}
	v_j & = & u   \hspace{.4in} \text {on} \hspace{.1in} \partial B_{2^{-j}}(z). 	
\end{array}\right.
$$

\vspace{.05in}
\noindent\\
By the same argument as before we can show $$ 	|Dv_k(z) -Du(z)| \leq C|z| \omega ( 2^{-k} ) \nonumber   $$

\vspace{.05in}
\noindent
Because $\Delta (u_k - v_k) = f(0)-f(z)$ in $B_{2^{-k}}(0) \cap B_{2^{-k}}(z)$ and $B_{2^{-k-1}}(z) \subset B_{2^{-k}}(0) \cap B_{2^{-k}}(z)$,
\vspace{.05in}
\begin{eqnarray*}
	|Dv_k(z) -Du_k(z)| & \leq & \| D(v_k - u_k)   \|_{L^{\infty} \left (  B_{2^{-k-2}}(z) \right ) }  \\
	& \leq & C \left ( 2^{k+1}   \| v_k - u_k   \|_{L^{\infty} \left (  B_{2^{-k-1}}(z) \right ) } + 2^{-k-1} |f(0)-f(z)|   \right )   \hspace{.2in} \big (\text{ by (\ref{eqn:estimate_Delta_with_constant}) \big )}\\
	& = & C \left ( 2^{k+1} \right )   \| v_k - u_k   \|_{L^{\infty} \left (  B_{2^{-k-1}}(z) \right ) } + C \left ( 2^{-k-1} \right )\omega(2^{-k-3}).
\end{eqnarray*}

\vspace{.05in}
\noindent
Then 
\begin{eqnarray}
	\label{eqn:u-minus_u_k_at_z_first_step}
	\left \vert   Du(z) - Du_k(z)  \right \vert & \leq & 	\left \vert   Du_k(z) - Dv_k(z)  \right \vert + 	\left \vert   Dv_k(z) - Du(z)  \right \vert  \nonumber \\
	& \leq & C \left ( 2^{k+1} \right )   \| v_k - u_k   \|_{L^{\infty} \left (  B_{2^{-k-1}}(z) \right ) } + C \left ( 2^{-k-1} \right )\omega(2^{-k-3}) + C|z| \omega ( 2^{-k} )  \nonumber  \\
	& \leq & C \left ( 2^{k+1} \right )   \| v_k - u_k   \|_{L^{\infty} \left (  B_{2^{-k-1}}(z) \right ) }  + C|z| \omega ( 2^{-k} )  .
\end{eqnarray}

\vspace{.05in}
\noindent
By (\ref{eqn:u-k_minus_u}) we know
$$ \| u_k - u   \|_{L^{\infty} \left (  B_{2^{-k-1}}(z) \right ) }  \leq   \| u_k - u   \|_{L^{\infty} \left (  B_{2^{-k}}(0) \right ) }  \leq  C \left ( 2^{-2k} \right ) \omega (2^{-k}),$$
and similarly we can prove  $$ \| v_k - u   \|_{L^{\infty} \left (  B_{2^{-k-1}}(z) \right ) }   \leq  C \left ( 2^{-2k} \right ) \omega (2^{-k}),$$ so

\begin{eqnarray*}
	\| u_k - v_k   \|_{L^{\infty} \left (  B_{2^{-k-1}}(z) \right ) } & \leq &  \| u_k - u   \|_{L^{\infty} \left (  B_{2^{-k-1}}(z) \right ) }  + \| v_k - u   \|_{L^{\infty} \left (  B_{2^{-k-1}}(z) \right ) } \\
	& \leq &  C \left ( 2^{-2k} \right ) \omega (2^{-k}). 	
\end{eqnarray*}

\vspace{.05in}
\noindent
Using this in (\ref{eqn:u-minus_u_k_at_z_first_step}) we have
\begin{eqnarray}
	\label{eqn:u-minus_u_k_at_z}
	\left \vert   Du(z) - Du_k(z)  \right \vert & \leq & 	C \left ( 2^{k+1} \right )  \left ( 2^{-2k} \right ) \omega (2^{-k})  +  C|z| \omega ( 2^{-k} )  \\
	& \leq &   C|z| \omega ( 2^{-k} ) . \nonumber
\end{eqnarray}

\vspace{.05in}
\noindent
Now we only need to estimate  $ |Du_k(z) - Du_k(0)|   $.  Let $$ h_j=u_j - u_{j-1}  \hspace{.4in} \text{for} \hspace{.1in} j=1,.., k.$$

\vspace{.05in}
\noindent
$h_j$ is harmonic, so by (\ref{eqn:estimate_Laplace})  $$ \| D^2 h_j \|_{L^{\infty} \left (  B_{2^{-j-1}}(0) \right )} \leq C \left ( 2^{2j}\right )  \|  h_j \|_{L^{\infty} \left (  B_{2^{-j}}(0) \right )}. $$

\vspace{.05in}
\noindent
Thus,
\allowdisplaybreaks
\begin{eqnarray*}
	\frac{|Dh_j(z) - D h_j(0)  |}{|z|} & \leq & \| D^2 h_j \|_{L^{\infty} \left (  B_{2^{-k-3}}(0) \right )}	\\
	& \leq & \| D^2 h_j \|_{L^{\infty} \left (  B_{2^{-j-1}}(0) \right )} \\
	& \leq & C \left ( 2^{2j}\right )  \|  h_j \|_{L^{\infty} \left (  B_{2^{-j}}(0) \right )} \\
	& = &  C \left ( 2^{2j}\right )  \|  u_j - u_{j-1} \|_{L^{\infty} \left (  B_{2^{-j}}(0) \right )} \\
	& \leq & C \left ( 2^{2j}\right )   \left (2^{-2(j-1)} \right ) \omega (2^{-(j-1)})  \hspace{.4in} \text{ by (\ref{eqn:u_K+1_minus_u_k}) }  \\
	& \leq & C\omega (2^{-(j-1)}).
\end{eqnarray*}

\vspace{.05in}
\noindent
Consequently, 
\allowdisplaybreaks
\begin{eqnarray*}
	|Du_k(z)-Du_k(0)| & \leq & |Du_1(z)-Du_1(0)| + \sum _{j=2}^k | Dh_j(z) - Dh_j(0)   |  \\
	& \leq & |Du_1(z)-Du_1(0)| + \sum _{j=2}^k C|z|\omega (2^{-j+1})\\
	& \leq & |z| \| D^2u_1  \|_{L^{\infty}(B_{ \frac{1}{4}  })} + C|z|\sum _{j=2}^k \omega (2^{-j+1}) .
\end{eqnarray*}

\vspace{.05in}
\noindent
Now we need to estimate $ \| D^2u_1  \|_{L^{\infty}(B_{\frac{1}{4}})} $. 

\vspace{.05in}
\noindent
Define a function $$\displaystyle \zeta (x) = u_1(x) - \frac{f(0)}{2n}|x|^2 + \frac{f(0)}{8n}.$$  

\vspace{.05in}
\noindent
Then $\zeta$ is harmonic because $\Delta \zeta =  \Delta u_1 - f(0) =0$, and $\zeta = u_1=u$ on $\partial B_{\frac{1}{2}}(0) $. 

\vspace{.05in}
\noindent
Furthermore, $D_{ij} \zeta = D_{ij} u_1$ when $i \neq j$, and $D_{ii} \zeta = D_{ii} u_1 - \frac{f(0)}{n}$.

\vspace{.05in}
\noindent
Therefore, 
\allowdisplaybreaks
\begin{eqnarray*}
	\| D^2u_1  \|_{L^{\infty}(B_{\frac{1}{4}})} & \leq &  \| D^2 \zeta \|_{L^{\infty}(B_{\frac{1}{4}})} + |f(0)| \\
	& \leq & C \|  \zeta \|_{L^{\infty}(B_{\frac{1}{2}})} + |f(0)| \\
	& = & C \|  \zeta \|_{L^{\infty}(\partial B_{\frac{1}{2}})} + |f(0)| \\
	& = & C \|  u \|_{L^{\infty}(\partial B_{\frac{1}{2}})} + |f(0)| \\
	& \leq & C \|  u \|_{L^{\infty}( B_{1})}+ |f(0)|.
\end{eqnarray*}

\vspace{.05in}
\noindent
It follows that 
\begin{equation}
	\label{eqn:u_k_z_minus_u_k_0}
	|Du_k(z) - Du_k(0)| \leq C|z|\|  u \|_{L^{\infty}( B_{1})} + |z||f(0)|+ C|z|\sum _{j=2}^k \omega (2^{-j+1}) .
\end{equation}

\vspace{.05in}
\noindent
Combining (\ref{eqn: 3_terms}), (\ref{eqn:estimate_u_k-u_at_0}), (\ref{eqn:u-minus_u_k_at_z}), and (\ref{eqn:u_k_z_minus_u_k_0}), we have
\begin{eqnarray*}
	|Du(z) - Du(0)| & \leq &  C|z| \omega ( 2^{-k} )  +  C|z| \|  u \|_{L^{\infty}( B_{1}) } + |z||f(0)|+ C|z|\sum _{j=2}^k \omega (2^{-j+1})\\
	& = & |z| \left (  \|  u \|_{L^{\infty}( B_{1}) } + |f(0)| +    \sum _{j=2}^{k+1} \omega (2^{-j+1}) \right ).
\end{eqnarray*}

\vspace{.05in}
\noindent
Finally, note that since $\omega(r)$ is increasing with $r$ increasing,
\begin{eqnarray*}
	\int_{|z|}^1 \frac{\omega(r)}{r} dr & \geq & \int_{\frac{1}{2^{k+3}}}^1 \frac{\omega(r)}{r} dr \\
	& \geq & \int_{\frac{1}{2^{k+3}}}^{\frac{1}{2^{k+2}}} \frac{\omega({\frac{1}{2^{k+3}}})}{r} dr + \int_{\frac{1}{2^{k+2}}}^{\frac{1}{2^{k+1}}} \frac{\omega({\frac{1}{2^{k+2}}})}{r} dr + \cdots + \int_{\frac{1}{2}}^1 \frac{\omega ({\frac{1}{2}}) }{r} dr  \\
	& = & (\ln 2) \omega \left ({\frac{1}{2^{k+3}}} \right ) + (\ln 2) \omega \left ({\frac{1}{2^{k+2}}} \right ) + \cdots + (\ln 2) \omega \left ({\frac{1}{2}} \right ).
\end{eqnarray*}

\noindent
Thus $$ \sum _{j=2}^{k+1} \omega (2^{-j+1}) <  \omega \left (\frac{1}{2} \right ) + \cdots +  \omega \left ({\frac{1}{2^{k+3}}} \right ) \leq C \int_{|z|}^1 \frac{\omega(r)}{r} dr. $$

\vspace{.05in}
\noindent
Therefore, we have proved that $$ 	|Du(z) - Du(0)|  \leq  |z| \left (  \|  u \|_{L^{\infty}( B_{1}) } + |f(0)| + \int_{|z|}^1 \frac{\omega(r)}{r} dr  \right ),$$ and this implies (\ref{eqn: first_derivative}).

\stop

\vspace{.4in}

\bibliographystyle{plain}
\bibliography{thesis}

\end{document}